\documentclass[12pt]{article}
\usepackage{amssymb,amsmath,amsfonts}
\usepackage{relsize}
\usepackage[sort]{cite}
\usepackage{esvect}
\usepackage{txfonts}
\usepackage{graphicx}
\newcommand{\ver}{{\rm ver}}
\newcommand{\vo}{{\rm vol}}

\newcommand{\mes}{{\rm mes_n}}

\newcommand{\dist}{{\rm dist}}
\newcommand{\simp}{{\rm simp_n}}

\newcommand{\conv}{{\rm conv}}

\newcommand{\vol}{{\rm vol}}

\newtheorem{theorem}{Theorem}[section]
\newtheorem{definition}[theorem]{Definition}
\newtheorem{corollary}[theorem]{Corollary}
\newtheorem{lemma}[theorem]{Lemma}
\newtheorem{remark}[theorem]{Remark}
\newtheorem{main problem}[theorem]{Main Problem}

\newcommand {\ip}[1]   {\langle{#1}\rangle}
\newcommand {\R}       {{\mathbb R}}
\newcommand {\RN}      {\R^n}
\newcommand {\smsk}    {\smallskip}
\newcommand {\msk}     {\medskip}
\newcommand {\bsk}     {\bigskip}

\newcommand {\rbx}     {\hspace{10mm}$\vartriangleleft$}

\newcommand {\SECT}[2] {\section*{\centerline{\normalsize
{\bf #1}}} \setcounter{section}{#2}
\setcounter{theorem}{0}\setcounter{equation}{0}}

\begin{document}
\parindent 1em
\parskip 0mm
\medskip
\centerline{\LARGE Optimal Lagrange Interpolation Projectors}
\vspace*{3mm}
\centerline{\LARGE and Legendre Polynomials}\vspace*{10mm}
\centerline{\large Mikhail Nevskii}\vspace*{5 mm}
\centerline {\it Department of Mathematics,  P.\,G.~Demidov Yaroslavl State University,}
\centerline {\it Sovetskaya str., 14, Yaroslavl, 150003, Russia}\vspace*{5 mm}
\centerline{\large May 2, 2024}

\vspace*{15mm}
\renewcommand{\thefootnote}{ }
\footnotetext[1]{\hspace{-6mm}
{\it E-mail address:} mnevsk55@yandex.ru}
\hrule\msk

\bsk
\par\noindent {\bf Abstract}
\bsk
\par Let $K$ be a convex body in $\RN$, and let  $\Pi_1(\RN)$ be the space of polynomials in $n$ variables
of degree at most $1$. Given an $(n+1)$-element set $Y\subset K$ in general position, we let $P_Y$ denote the Lagrange interpolation projector $P_Y: C(K)\to \Pi_1(\RN)$ with nodes in $Y$. In this paper, we study upper and lower bounds for the norm of the optimal Lagrange interpolation projector, i.e., the projector with minimal operator norm where the minimum is taken over all $(n+1)$-element sets of interpolation nodes in $K$. We denote this minimal norm by $\theta_n(K)$.
\par Our main result, Theorem \ref{th_theta_n_K_thru_vol_K_simp_K},
 provides an explicit lower bound for the constant $\theta_n(K)$ for an arbitrary convex body $K\subset\RN$ and an arbitrary $n\ge 1$. We prove that
$\theta_n(K)\ge \chi_n^{-1}\left({\vo(K)}/{\simp(K)}\right)$ where $\chi_n$ is the Legendre polynomial of degree $n$ and
$\simp(K)$ is the maximum volume of a simplex contained in $K$. The proof of this result relies on a geometric
characterization of the Legendre polynomials in terms of  the volumes of certain convex polyhedra. More specifically,
we show that for every $\gamma\ge 1$ the volume of the set
$\left\{x=(x_1,...,x_n)\in{\mathbb R}^n :
\sum |x_j| +\left|1- \sum x_j\right|\le\gamma\right\}$ is equal to ${\chi_n(\gamma)}/{n!}$.
\par Furthermore, if $K$ is an $n$-dimensional ball, this approach leads us to the
equivalence
  $\theta_n(K) \asymp\sqrt{n}$ 
 which is complemented by the exact formula  for $\theta_n(K)$. If $K$ is an $n$-dimensional cube, we obtain
explicit efficient formulae for upper and lower bounds of the constant $\theta_n(K)$; moreover, for small $n$, these estimates enable us to compute the exact values of this constant.
\bsk
{\small
\smsk
}
\msk
\hrule

\bsk


\SECT{1. Introduction}{1}\label{intro}

\par Let $K$ be a convex body in $\RN$.
 Given a function $f\in C(K)$, we let $E_1(f;K)_{C(K)}$ denote the best approximation of $f$ (in the $C(K)$-norm) by polynomials of degree at most $1$ defined on $\RN$. (We~denote the~space of these polynomials by $\Pi_1(\RN)$.) 
\par Let $P$ be a polynomial projector on $K$, i.e., a linear operator from $C(K)$ into $\Pi_1(\RN)$ such that
$P(Pf|_K)=Pf$ for every function $f\in C(K)$. Then the following well known inequality
\begin{equation}\label{Lebesgue_Lemma}
E_1(f;K)_{C(K)}\le \|f-Pf\|_{C(K)}\le (1+\|P\|_K)E_1(f;K)_{C(K)}.
\end{equation}
holds. Here $\|P\|_K=\sup\{\|Pf|_K\|_{C(K)}:\|f\|_{C(K)}\le 1\}$ is the $C(K)$ operator norm of $P$.
\smsk
\par Clearly, the first inequality in (\ref{Lebesgue_Lemma}) is trivial because $Pf\in\Pi_1(\RN)$. The second inequality is known in the  literature as Lebesgue's Lemma. See, e.g., \cite[p.30]{devore_lorentz_1993}. This lemma shows that, if the norm $\|P\|_K$ is not large, there is no significant loss of  accuracy if we replace the best $C(K)$-approximant of $f$ by~the~value $Pf$ {\it which linearly depends on the function $f$}.
\smsk
\par This approach provides a (fairly simple) tool for ``linearizing'' the best approximation element in~various problems of Analysis involving the approximation of continuous functions defined on subsets of $\RN$. It also motivates us to study the following problem.
\begin{main problem}\label{MP} {\em Let $K$ be a convex body in $\RN$, and let $Y$ be an $(n+1)$-element subset of $K$ in~``general position'' (i.e., $Y$ is the family of vertices of a nondegenerate simplex in $\RN$). Let $P_Y$ be the first order Lagrange interpolation projector with nodes in $Y$ (i.e., $P_Y(f)\in\Pi_1(\RN)$ and $P_Y(f)=f$ on $Y$ for every $f\in C(K)$).
\smsk
\par (i) Find an $(n+1)$-element set $Y_{opt}\subset K$ such that $P_{Y_{opt}}$ has the smallest operator norm among all Lagrange interpolation projectors $P_Y$ where $Y$ is an arbitrary $(n+1)$-element subset of $K$.
\par We refer to the Lagrange interpolation projector $P_{Y_{opt}}$ as the optimal interpolation projector, and we call the set $Y_{opt}$ as an optimal set of interpolation nodes.
\smsk
\par (ii) Calculate or give efficient estimates for the upper and lower bounds of the constant $$\theta_n(K)=\|P_{Y_{opt}}\|_K,$$ i.e., the minimum of the operator norm $\|P_Y\|_K$ where the minimum is taken over all $(n+1)$-element sets $Y$ of interpolation nodes in $K$.}
\end{main problem}
\par Let us note that both upper and lower bounds of the operator norms of interpolation projectors are important in applications. Usually, upper estimates of minimal norms of projectors are obtained by~considering projectors of some special type. The technique of obtaining lower estimates of minimal norms is fundamentally different -- in this case it is necessary to prove the lower estimate for an~arbitrary projector.
\par The construction and evaluation of interpolation projectors is a classical topic in Approximation Theory. These problems have been treated in many papers and monographs, see, e.g., \cite{barthelmann_2000,deboor_1994,demarchi_2015,
gasca_2000,gunzburger_2014,kaarnioja_2015,
pashkovskij_1983,rivlin_1974}.
\smsk
\par In this paper, we show that there exists a Lagrange interpolation projector whose operator norm does not exceed $n+1$ (Theorem \ref{th_norm_P_for_max_simplex}). On the other hand, we prove that for any Lagrange interpolation projector $P:C(K)\to \Pi_1(\RN)$ the following inequality
\begin{equation}\label{LEG-POL}
\|P\|\geq
\chi_n^{-1}
\left(\frac{\vo(K)}{\simp(K)}\right)
\end{equation}
holds. Here $\chi_n$ is the standardized Legendre polynomial of degree $n$ and  $\simp(K)$ is the maximum volume of a simplex contained in~$K$. See  Theorem \ref{th_theta_n_K_thru_vol_K_simp_K}.
\par The main point of the proof of inequality (\ref{LEG-POL}) is the following geometrical characterization of Legendre polynomials given in Theorem \ref{th_mes_E_xi_eqs}. This theorem states that for $\gamma\ge 1$ the volume of the set
$$
\Bigl\{ x\in{\mathbb R}^n :
\sum
_{j=1}^n
|x_j| +
\Bigl|1- \sum
_{j=1}^n
x_j\Bigr| \leq \gamma \Bigr\}
$$
is equal to ${\chi_n(\gamma)}/{n!}$.
\smsk
\par In the cases when $K$ coincides with an $n$-dimensional cube or ball, the above inequality yields
the~estimate $\|P\|\geq c\sqrt{n}$.
\smsk
\par Let us describe the content of the paper in more detail. In Section 2 we give basic notation, definitions, and  preliminary information. Section 3 contains upper estimates of the minimum absorption coefficient of  a convex body $K$  by a simplex and also of the minimum projector norm for linear interpolation on $K$. In Section 4 we prove Theorem \ref{th_mes_E_xi_eqs} mentioned above. In Section 5 we prove inequality (\ref{LEG-POL}). In Section 6 we present several explicit lower bounds for the constant $\theta_n(K)$ provided $K$ is an~$n$-dimensional ball or an $n$-dimensional cube. Section 7 contains inequalities of the considered type for interpolation by linear functions on an arbitrary compact set in $\RN$. Finally, Section 8 contains concluding remarks, a review of some results on this topic, and a description of open questions.

\SECT{2. Notation and preliminaries}{2}
\label{main_defs_and_denots}


\indent\par Let us fix some notation. Throughout the paper  $n\in{\mathbb N}$ is a positive integer.
 Given $x=(x_1,...x_n)\in~\RN$, by $\|x\|$ we denote its standard Euclidean norm
$$
\|x\|=\sqrt{\langle x,x\rangle}=\left(\sum\limits_{i=1}^n x_i^2\right)^{1/2}.
$$
Hereafter, for $x=(x_1,...x_n), y=(y_1,...y_n)\in\RN$ by
$
\ip{x,y}
$
we denote the standard inner product in $\RN$:
$
\ip{x,y}=x_1y_1+...+x_ny_n
$

\par Given $x^{(0)}\in\RN$ and $R>0$, we let $B(x^{(0)};R)$ denote the Euclidean ball with center $x^{(0)}$ and radius $R$:
$$
B\left(x^{(0)};R\right)=\{x\in{\mathbb R}^n: \|x-x^{(0)}\|\leq R \}.
$$

We also set
$$
B_n=B(0;1)~~~\text{and}~~~
Q_n=[0,1]^n.
$$


The notation $L(n)\asymp M(n)$ means that there exist absolute constants $c_1,c_2>0$ such that
\linebreak $c_1M(n)\leq L(n)\leq c_2 M(n)$.

Let $K$ be {\it a convex body in ${\mathbb R}^n$}, i.e., a compact convex subset of ${\mathbb R}^n$
with nonempty interior.
The~symbol $\sigma K$  denotes a homothetic copy of $K$ with  homothety  center  at
the center of gravity of $K$
and homothety coefficient  $\sigma.$

We let $\vol(K)$ denote the volume of $K$.
If $K$ is a convex polytope, then by $\ver(K)$ we denote the set of all vertices of $K$. By {\it a translate} we mean the result of a parallel shift.

We say that an $n$-dimensional simplex $S$ is {\it circumscribed around a convex body $K$} if $K\subset S$ and~each $(n-1)$-dimensional face of  $S$ contains a point of $K$. A convex polytope is {\it inscribed into $K$} if every vertex of this polytope belongs to the boundary of
$K$.

\begin{definition}\label{def_ax_diam}{\em
Let $i\in\{1,...,n\}$ and let $d_i(K)$ be the maximal length of a segment contained in $K$ and parallel to the
$x_i$-axis. We refer to $d_i(K)$  as  the $i$th
axial diameter of $K$.}
\end{definition}

\smallskip
The notion of {\it axial diameter of a convex body} was introduced by Scott \cite{scott_1985, scott_1989}.

\begin{definition}\label{def_xi_alpha_K_1_K_2}{\em Given convex bodies $K_1$, $K_2$, by $\xi(K_1;K_2)$ we denote the minimal $\sigma\geq 1$ with the~pro\-perty
$K_1\subset \sigma K_2$.  We call $\xi(K_1,K_2)$ {\it the absorption index} of $K_1$ by $K_2$.

\par By $\alpha(K_1,K_2)$ we denote the minimal $\sigma>0$ such that $K_1$ is a subset of a translate of $\sigma K_2$. }
\end{definition}

\smallskip
\par Note that the equality
$\xi(K_1;K_2)$ $=1$ is equivalent to the inclusion  $K_1\subset K_2$. Clearly, $\alpha(K_1,K_2)$ $\leq $ $\xi(K_1,K_2)$.

\par\begin{definition}\label{xi_n_K_def}{\em By $\xi_n(K)$ we denote the minimal absorption index of a convex body $K$ by an inner nondegenerate simplex. In other words,
$$
\xi_n(K)=\min \{ \xi(K;S): \,
S  \mbox{ is an $n$-dimensional simplex,} \,
S\subset K, \, \vo(S)\ne 0\}.
$$
}
\end{definition}

\smallskip
By $C(K)$ we denote the space of all continuous functions
$f:K\to{\mathbb R}$ with the uniform
norm
$$
\|f\|_{C(K)}=\max\limits_{x\in K}|f(x)|.
$$
We let $\Pi_1\left({\mathbb R}^n\right)$ denote the space of polynomials in $n$ variables
of degree at most $1$.

Let $S$ be a nondegenerate simplex in ${\mathbb R}^n$ with vertices $x^{(j)}=\left(x_1^{(j)},\ldots,x_n^{(j)}\right),$
$1\leq j\leq n+1.$ We define {\it the vertex matrix} of this simplex by
$$
{\bf A}=
\left( \begin{array}{cccc}
x_1^{(1)}&\ldots&x_n^{(1)}&1\\
x_1^{(2)}&\ldots&x_n^{(2)}&1\\
\vdots&\vdots&\vdots&\vdots\\
x_1^{(n+1)}&\ldots&x_n^{(n+1)}&1\\
\end{array}
\right).
$$

Clearly, matrix ${\bf A}$ is nondegenerate and
\begin{equation}\label{vol_S_det_A_eq}
\vol(S)=\frac{|\det({\bf A})|}{n!}.
\end{equation}
Let
\begin{equation}\label{ACV}
{\bf A}^{-1}=\left(l_{ij}\right)_{i,j=1}^{n+1}.
\end{equation}

\begin{definition}\label{basic_Lagr_polynomials}{\em Linear  polynomials
$$
\lambda_j(x)=
l_{1j}x_1+\ldots+
l_{nj}x_n+l_{n+1,j},~~~~j=1,...,n+1,
$$
are called {\it the~basic Lagrange polynomials} corresponding to $S$.}
\end{definition}

\smallskip
These polynomials have the following property:
$$
\lambda_j\left(x^{(k)}\right)=\delta_j^k~~~\text{for all}~~~j,k=1,...,n+1.
$$
Here $\delta_j^k$ is the Kronecker delta.
For an arbitrary $x\in{\mathbb R}^n$, we have
$$
x=\sum_{j=1}^{n+1} \lambda_j(x)x^{(j)}, \quad \sum_{j=1}^{n+1} \lambda_j(x)=1.
$$
Thus, $\lambda_j(x)$ are {\it the barycentric coordinates of $x$} with respect to the simplex $S$.
In turn, equations $\lambda_j(x)=0$, $j=1,...,n+1$, define the $(n-1)$-dimensional hyperplanes containing the faces of $S.$  Therefore,
$$
S=\left \{ x\in {\mathbb R}^n: \, \lambda_j(x) \geq 0, \, j=1,\ldots,n+1
\right\}.
$$

Also let us note that for every $j=1,...,n+1$ we have
\begin{equation}\label{Lagr_pol_thru_dets}
\lambda_j(x)=\frac{\Delta_j(x)}{\Delta}.
\end{equation}
Here $\Delta=\det({\bf A})$ and $\Delta_j(x)$ is the determinant that appears from $\Delta$ after replacing the $j$th row by~the~row $(x_1\  \ldots\  x_n \ 1).$
For more information on $\lambda_j$, see ~\cite{nevskii_monograph}, \cite{nev_ukh_posobie_2020}.

In \cite{nevskii_dcg_2011} we show that
\begin{equation}\label{xi_K_S_formula}
\xi(K;S)=(n+1)\max_{1\leq k\leq n+1}
\max_{x\in K}(-\lambda_k(x))+1 \quad (K\not\subset S),
\end{equation}
\begin{equation}\label{alpha_K_S_general_formula}
\alpha(K;S)=\sum_{j=1}^{n+1} \max_{x\in K} (-\lambda_j(x))+1
\end{equation}
(see also \cite{nevskii_monograph}). The equality $\xi(K;S)=\alpha(K;S)$  holds true if and only if the simplex $\xi(S)S$ is circumscribed around $K$. This  is also equivalent to the relation
\begin{equation}\label{xi_S_S_circ_condition_around_K}
 \max_{x\in K} \left(-\lambda_1(x)\right)=
\ldots=
\max_{x\in K} \left(-\lambda_{n+1}(x)\right).
\end{equation}

If $K$ is a convex polytope, then the maxima on $K$ in
\eqref{xi_K_S_formula}--\eqref{xi_S_S_circ_condition_around_K} can also be taken over $x\in \ver(K)$. Note that for $K=Q_n$
formula \eqref{xi_K_S_formula} is proved in \cite{nevskii_matzam_2010}.

 For an arbitrary $n$-dimensional nondegenerate simplex $S$,
\begin{equation}\label{alpha_d_i_S_eq}
\alpha(Q_n;S)=\sum_{i=1}^n\frac{1}{d_i(S)}
\end{equation}
Two various approaches to \eqref{alpha_d_i_S_eq} are given in \cite{nevskii_dcg_2011}.

\smallskip
\begin{definition}\label{interp_projector_def}{\em
Let $x^{(j)}\in K$, $1\leq j\leq n+1,$ be  the vertices of a nondegenerate simplex $S$. The interpolation projector $P:C(K)\to \Pi_1({\mathbb R}^n)$  with nodes $x^{(j)}$ is determined by the equalities $Pf\left(x^{(j)}\right)=f_j=f\left(x^{(j)}\right)$ , $1\leq j\leq n+1.$
We say that an interpolation projector $P:C(K)\to \Pi_1({\mathbb R}^n)$ and a nondegenerate
simplex $S\subset K$ {\it correspond to each other} if
the  nodes of $P$ coincide with  the~vertices of $S$. We use notation $P_S$ and $S_P$.}
\end{definition}

\smallskip
For an interpolation projector $P=P_S$, the analogue of Lagrange interpolation formula holds:
\begin{equation}\label{interp_Lagrange_formula}
Pf(x)=\sum\limits_{j=1}^{n+1}
f\left(x^{(j)}\right)\lambda_j(x).
\end{equation}
Here $\lambda_j$ are the basic Lagrange polynomials of the simplex $S_P$ (see Definition 1.4).

We let $\|P\|_K$ denote the norm of $P$ as an operator from $C(K)$ into $C(K)$.  Thanks to (\ref{interp_Lagrange_formula}),
\begin{eqnarray}
\|P\|_K &=& \sup_{\|f\|_{C(K)}=1} \|Pf\|_{C(K)}=
\sup_{-1\leq f_j\leq 1} \max_{x\in K}\left| \sum_{j=1}^{n+1} f_j\lambda_j(x)\right| \nonumber\\
&=&
\max_{x\in K}\sup_{-1\leq f_j\leq 1}\left| \sum_{j=1}^{n+1} f_j\lambda_j(x)\right|=
\max_{x\in K}\sup_{-1\leq f_j\leq 1}\sum_{j=1}^{n+1} f_j\lambda_j(x). \nonumber
\end{eqnarray}
Because $\sum f_j\lambda_j(x)$ is linear in $x$ and
$f_1,\ldots,f_{n+1}$, we have
\begin{equation}\label{norm_P_intro_cepochka}
\|P\|_K=
\max_{x\in K}\max_{f_j=\pm 1} \sum_{j=1}^{n+1}
f_j\lambda_j(x)
=\max_{x\in K}\sum_{j=1}^{n+1}
|\lambda_j(x)|.
\end{equation}
If $K$ is a convex polytope in ${\mathbb R}^n$
(e.\,g., $K$ is a cube), a simpler equality
\begin{equation}\label{norm_P_cube_formula}
\|P\|_K= \max_{x\in\ver(K)}\sum_{j=1}^{n+1}
|\lambda_j(x)|
\end{equation}
holds.

\begin{definition}\label{theta_n_K_def}{\em We let $\theta_n(K)$ denote the minimal value of $\|P_S\|_K$ where $S$ runs over all nondegenerate simplices with vertices in $K$. An interpolation projector $P:C(K)\to\Pi_1\left({\mathbb R}^n\right)$ is called minimal if $\|P\|_K=\theta_n(K)$.}
\end{definition}

\smsk
\smallskip
It was shown in \cite{nevskii_mais_2008_15_3} that for any interpolation projector
$P:C(K)\to\Pi_1\left({\mathbb R}^n\right)$ and the corresponding simplex  $S$ we have
\begin{equation}\label{nev_ksi_P_ineq}
\frac{n+1}{2n}\Bigl( \|P\|_K-1\Bigr)+1\leq
\xi(K;S) \leq
\frac{n+1}{2}\Bigl( \|P\|_K-1\Bigr)+1.
\end{equation}

\noindent Thanks to \eqref{nev_ksi_P_ineq},
\begin{equation}\label{nev_ksi_n_K_theta_n_K_ineq}
\frac{n+1}{2n}\Bigl( \theta_n(K)-1\Bigr)+1\leq
\xi_n(K) \leq
\frac{n+1}{2}\Bigl( \theta_n(K)-1\Bigr)+1.
\end{equation}
Obviously, if a projector $P$ satisfies the equality
\begin{equation}\label{ksi_n_K_norm_P_eq}
\xi_n(K) =
\frac{n+1}{2}\Bigl( \|P\|_K-1\Bigr)+1,
\end{equation}
then $P$ is minimal and the right-hand relation in \eqref{nev_ksi_n_K_theta_n_K_ineq} becomes an equality.

Occasionally, we will consider the case when $n+1$ is~an~Hadamard number, i.e., there exists an~Hadamard matrix of~order $n+1$.

\begin{definition}\label{def_Had_matr}{\em An  Hadamard matrix of order $m$ is a square binary matrix $\bf H$ with entries either $1$ or $-1$ satisfying the  equality
$${\bf H}^{-1}=\frac{1}{m}\, {\bf H}^{T}.$$
An integer $m$, for which an Hadamard matrix of order $m$ exists, is called an Hadamard number.}
\end{definition}

\smallskip
Thus, the rows of ${\bf H}$ are pairwise orthogonal with respect to the standard scalar product on ${\mathbb R}^m.$

The order  of an Hadamard matrix is $1$ or $2$ or some multiple of $4$ (see \cite{hall_1970}).  But it is still unknown whether an Hadamard matrix exists for every order of the form $m=4k$. This  is one of the longest lasting open problems in Mathe\-matics called {\it the Hadamard matrix conjecture}.
The orders below $1500$ for which Hadamard matrices are not known are $668, 716, 892, 956, 1132, 1244, 1388$, and $1436$ (see, e.g., \cite{horadam_2007, manjhi_2022}).

\begin{definition}\label{h_n_nu_n_def}{\em
 By $h_n$ we denote the maximum value of a determinant of order $n$ with entries $0$ or~$1$. Let $\nu_n$ be the maximum volume of an $n$-dimensional simplex
 contained in $Q_n$.}
 \end{definition}

 \smallskip
 The numbers $h_n$ and $\nu_n$ satisfy  the equality $h_n=n!\nu_n$ (see \cite{hudelson_1996}).
For any $n$, there exists  in $Q_n$ a~maximum volume simplex with some vertex coinciding with a vertex of the cube.
For $n>1$, the following inequalities hold
\begin{equation}\label{adamar_clements_lindstrem1}\frac{1}{2}\left( 1-\frac{\log(4/3)}{\log n}\right ) n \log n
< \log(2^{n-1}h_{n-1})\leq \frac{1}{2} \, n\log n.
\end{equation}
The right-hand inequality in
(\ref{adamar_clements_lindstrem1}) was proved by Hadamard
\cite{hadamard_1893} and the left-hand one by Clements and  Lindstr\"om \cite{clements_lindstrem_1965}.
Consequently, for all $n\in {\mathbb N}$
\begin{equation}\label{adamar_clements_lindstrem2}
\left(\frac{3}{4}\right)^{(n+1)/2}\,\frac{\left(n+1\right)^{(n+1)/2}}{2^n}<
h_n\leq\frac{\left(n+1\right)^{(n+1)/2}}{2^n},
\end{equation}
\begin{equation}\label{adamar_clements_lindstrem}
\left(\frac{3}{4}\right)^{(n+1)/2}\,\frac{\left(n+1\right)^{(n+1)/2}}{2^nn!}<
\nu_n\leq\frac{\left(n+1\right)^{(n+1)/2}}{2^nn!}.
\end{equation}

The right-hand equality in each relation
holds if and only if $n+1$ is an Hadamard number \cite{hudelson_1996}.
In some cases the right-hand inequality in
(\ref{adamar_clements_lindstrem})
has been improved.
For instance, if $n$ is even, then
\begin{equation}\label{nu_n_n_even}
\nu_n\leq \frac{n^{n/2}\sqrt{2n+1}}{2^nn!}.
\end{equation}

If $n>1$ and $n\equiv 1({\rm mod}~4),$ then
\begin{equation}\label{nu_n_n_odd}
\nu_n\leq \frac{(n-1)^{(n-1)/2}}{2^{n-1}(n-1)!}
\end{equation}
(see \cite{hudelson_1996}).
For many $n$, the  values of $\nu_n$ and $h_n$ are known exactly.  The first $12$ numbers $\nu_n$ are
$$
\nu_1=1, \quad \nu_2= \frac{1}{2}, \quad \nu_3= \frac{1}{3},
\quad \nu_4=\frac{1}{8}, \quad \nu_5=\frac{1}{24}, \quad
\nu_6=\frac{1}{80}, \quad \nu_7=\frac{2}{315},$$\
$$\nu_8=\frac{1}{720}, \quad \nu_9=\frac{1}{2520}, \quad
\nu_{10}=\frac{1}{11340}, \quad \nu_{11}=\frac{9}{246400},
\quad \nu_{12}=\frac{3}{394240}.
$$

\begin{definition}\label{kappa_n_sigma_n_def}{\em
Let $\varkappa_n$ be the volume of the ball $B_n$, and let $\sigma_n $ be the volume of a regular simplex inscribed into $B_n$.}
\end{definition}

The numbers $\varkappa_n$ and  $\sigma_n$ are known exactly. Namely,
\begin{equation}\label{kappa_n_sigma_n_whole_1}
\varkappa_n=\frac{\pi^
{{n}/{2}}}
{\Gamma\left({n}/{2}+1\right)},\qquad
\sigma_n=\frac{1}{n!}\sqrt{n+1}\left(\frac{n+1}{n}\right)^{{n}/{2}},
\end{equation}
\begin{equation}\label{kappa_n_even_and_odd_1}
\varkappa_{2k}=\frac{\pi^{k}}{k!},\qquad
\varkappa_{2k+1}=\frac{2^{k+1}\pi^{k}}{(2k+1)!!}=
\frac{2(k!)(4\pi)^k}{(2k+1)!}.
\end{equation}

\begin{definition}\label{simp_n_K_def}{\em
We let $\simp(K)$ denote the maximum volume of a nondegenerate simplex $S$ with vertices in $K$.}
\end{definition}

Obviously, $\simp(Q_n)=\nu_n$.  A regular simplex inscribed into an $n$-dimensional ball has the maximum volume among all simplices contained in this ball. There are no other simplices with this property  (see \cite{fejes_tot_1964, slepian_1969,
vandev_1992}). Thus, $\simp(B_n)=\sigma_n$.

\SECT{3. Inequalities~~ $\xi_n(K)\leq n+2$~ and~ $\theta_n(K)\leq n+1$}{3}\label{ksi_n_leq_nplus2_theta_n_K_leq_nplus1}

\begin{theorem}\label{th_norm_P_for_max_simplex}
Let $K$ be a convex body in $\RN$, and let $S$ be an arbitrary maximum volume simplex contained  in $K$.
Then
\begin{equation}\label{th_xi_S_norm_PS_ineqs_for_maxvol_S}
\xi(K;S)\leq n+2~~~~~\text{and}~~~~\quad \|P_S\|_K\leq n+1.
\end{equation}
\end{theorem}

\smallskip
{\it Proof.} Using a purely geometric approach,  M. Lassak  \cite{lassak_beitr_2011} have proved that for every maximum volume simplex $S$ in $K$ the following inclusions
\begin{equation}\label{S_sub_K_sub_n_plus_2_S}
S\subset K\subset (n+2)S
\end{equation}
hold. From this and Definition \ref{def_xi_alpha_K_1_K_2}, we have $\xi(K;S)\leq n+2$.

\par Let us note that \eqref{S_sub_K_sub_n_plus_2_S} also follows from formula \eqref{xi_K_S_formula}. In fact, because simplex $S\subset K$ has maximum volume, $|\Delta_j(x)|\leq |\Delta|$ for any $j=1,\ldots,n+1$ and $x\in K.$
(This inequality is immediate from \eqref{vol_S_det_A_eq}.) Thanks to \eqref{Lagr_pol_thru_dets},
\begin{equation}\label{LM-1}
-\lambda_j(x)\leq |\lambda(x)|=\frac{|\Delta_j(x)|}{|\Delta|} \leq 1, \quad x\in K.
\end{equation}
(Recall that $\lambda_j$ are the basic Lagrange polynomials for the simplex $S$. See Definition \ref{basic_Lagr_polynomials}.)
By \eqref{xi_K_S_formula},
$$\xi(K;S)=(n+1)\max_{1\leq k\leq n+1}
\max_{x\in K}(-\lambda_k(x))+1\leq n+2$$
proving the left-hand inequality in \eqref{th_xi_S_norm_PS_ineqs_for_maxvol_S}.

\par The right-hand inequality in \eqref{th_xi_S_norm_PS_ineqs_for_maxvol_S} follows from \eqref{norm_P_intro_cepochka}. Indeed, because $|\lambda(x)| \leq 1,$ we have
$$\|P_S\|_K
=\max_{x\in K}\sum_{j=1}^{n+1}
|\lambda_j(x)| \leq n+1.$$

\par The proof of the theorem is complete.\hfill$\Box$

\par The following corollary is immediate from Theorem \ref{th_xi_S_norm_PS_ineqs_for_maxvol_S}.
\begin{corollary}\label{corol_xi_n_K_leq_n+2_theta_n_K_leq_n+1}
For any covex body $K\subset \RN$, 
\begin{equation}\label{xi_n_K_leq_n+2_theta_n_K_leq_n+1}
\xi_n(K) \leq n+2,~~~ \quad \theta_n(K)\leq n+1.
\end{equation}
\end{corollary}

\SECT{4. Legendre polynomials and the measure of $E_{n,\gamma}$}{4}\label{legendre_pol_mes_E}

\indent\par {\it The standardized Legendre polynomial of degree $n$} is the function
$$
\chi_n(t)=\frac{1}{2^nn!}\left[ (t^2-1)^n \right] ^{(n)},~~~~~~t\in\R.
$$
(Rodrigues' formula). For the properties of $\chi_n$ see, e.g., \cite{suetin_1979,sege_1975}.
Legendre polynomials are orthogonal on the segment
$[-1,1]$
with the weight $w(t)\equiv1.$
The first  standardized  Legendre polynomials are
$$\chi_0(t)=1, \quad
\chi_1(t)=t, \quad
\chi_2(t)=\frac{1}{2}\left(3t^2-1\right), \quad
\chi_3(t)=\frac{1}{2}\left(5t^3-3t\right),$$
$$\chi_4(t)=\frac{1}{8}\left(35t^4-30t^2+3\right), \quad
\chi_5(t)=\frac{1}{8}\left(63t^5-70t^3+15t\right).$$

\par For these polynomials the following recurrent relations
\begin{equation}\label{reccurent_Legendre}
\chi_{n+1}(t)=
\tfrac{2n+1}{n+1}\,t\cdot\chi_n(t)
-\tfrac{n}{n+1}\,\chi_{n-1}(t)
\end{equation}
hold. In particular, this implies
$\chi_n(1)=1$. Let us also note that, if $n\ge 1$ then $\chi_n(t)$ increases as~$t\geq 1$. (This property also easily follows from formula \eqref{th_2_1_1} given below.)

\par We let $\chi_n^{-1}$ denote the function inverse to $\chi_n$ on the semi-axis $[1,+\infty)$.
\smsk

\par One of the key statements of our approach to the problems of optimal Lagrangian interpolation is Theorem \ref{th_mes_E_xi_eqs} formulated and proven in this section. This theorem reveals some rather surprising connections between Legendre polynomials with the problems of optimal Lagrangian interpolation.

\par Given $\gamma\geq 1$, we let $E_{n,\gamma}$ define a set
\begin{equation}\label{ENG}
E_{n,\gamma}= \left\{ x\in{\mathbb R}^n :
\sum_{j=1}^n |x_j| +
\left|1- \sum_{j=1}^n x_j\right| \leq \gamma \right\}.
\end{equation}

\smallskip

\begin{theorem}\label{th_mes_E_xi_eqs}
The  following equalities hold:
\begin{equation}\label{th_2_1_1}
\mes(E_{n,\gamma})
= \frac{1}{2^n n!} \sum_{i=0}^{n}  {n \choose i}^2 (\gamma - 1)^{n-i}
(\gamma + 1)^i = \frac{\chi_n(\gamma)}{n!}.
\end{equation}
\end{theorem}

\smallskip
\par This result was obtained in \cite{nevskii_mais_2003_10_1}. Unfortunately, this paper is practically inaccessible to a wide audience, so for the convenience of the reader we present it here.

\smallskip
{\it Proof.} First, let us  prove the left-hand equality in \eqref{th_2_1_1}.
Let
$$
E^{(1)}=\{x\in E_{n,\gamma}: \sum x_i>1 \}~~~~~
\text{and}~~~~
E^{(2)}=\{x\in E_{n,\gamma}:\sum x_i \leq 1\}.
$$

Let us give explicit formulae for the volumes   $m_1=\mes(E^{(1)})$ and $m_2=\mes(E^{(2)}).$

Let us fix $k$, $1\le k\le n$, and consider a non-empty subset $G\subset E^{(1)}$ consisting of all points $x=(x_1,...,x_n)$ such that
$x_1,\ldots,x_k\geq 0$ and $x_{k+1},\ldots, x_n< 0$. Let $y_i=x_i$ for $i=1,...,k$,
$y_{i}=-x_{i}$ for~$i=k+1,...,n$, and let $y=(y_1,...,y_n)$. Then
$$G=\left\{y: 1+y_{k+1}+\ldots+y_n\leq y_1+\ldots+y_k\leq \frac{\gamma+1}{2}, \
y_i\geq 0 \right \},$$
hence
\begin{eqnarray}
\mes(G) &=&
\int\limits_{1}^{\alpha} \,dy_1
\int\limits_{1}^{\alpha-y_1} \,dy_2
\ldots
\int\limits_{1}^{\alpha-y_1-\ldots-y_{k-1}} \,dy_k
\nonumber\\ \nonumber \\
& &
\int\limits_{0}^{y_1+\ldots+y_k-1} \,dy_{k+1}
\int\limits_{0}^{y_1+\ldots+y_k-1-y_{k+1}} \,dy_{k+2}
\ldots
\int\limits_{0}^{y_1+\ldots+y_k-1-y_{k+1}-\ldots-y_{n-1}} \,dy_n. \nonumber
\end{eqnarray}
Throughout the proof, $\alpha=
(\gamma+1)/2.$
If $b>0$, then
$$\int\limits_{0}^{b} \,dz_1
\int\limits_{0}^{b-z_1} \,dz_2\ldots\int\limits_{0}^{b-z_1-\ldots-z_{l-1} }
\,dz_l=\frac{b^l}{l!},$$
so
\begin{eqnarray}
\mes(G) &=&
\int\limits_{1}^{\alpha} \,dy_1
\int\limits_{1}^{\alpha-y_1} \,dy_2
\ldots
\int\limits_{1}^{\alpha-y_1-\ldots-y_{k-1}}
\frac{1}{(n-k)!}{\left(y_1+\ldots+y_k-1\right)}^{n-k}
\,dy_k \nonumber\\ \nonumber \\
&=& \left (\int\limits_{\ y_1+\ldots+y_k\leq\alpha} -\int\limits_{y_1+\ldots+y_k\leq 1}
\right)
\frac{1}{(n-k)!}  {(y_1+\ldots+y_k-1)}^{n-k}
\,dy_1\ldots \,dy_k \nonumber\\ \nonumber\\
&=&  J_1 - J_2.
\nonumber
\end{eqnarray}
The first integral equals
$$J_1=
\sum_{j=1}^{k} (-1)^{j+1} \frac{(\alpha-1)^{n-k+j}}{(n-k+j)!}
\frac{\alpha^{k-j}}{(k-j)!}+\frac{(-1)^{n+k}}{n!}
.$$
The value of $J_2$ appears from this expression
if instead of $\alpha$ we take 1. Consequently,

\begin{eqnarray}
\mes(G)&=& J_1 - J_2\nonumber\\ &=&
\sum_{j=1}^{k} (-1)^{j+1} \frac{ (\alpha-1)^{n-k+j} }{(n-k+j)!}
\frac{ \alpha^{k-j} }{ (k-j)! }
\nonumber\\
&=& \frac{ (-1)^{k+1} }{n!}
\sum_{i=0}^{k-1} {n\choose i}
(\alpha-1)^{n-i}(-\alpha)^i.
\end{eqnarray}

Clearly, the set $E^{(1)}$ is the union of all pairwise disjoint sets $G$ with various
$k=1,\ldots,n,$. Therefore, the measure of $E^{(1)}$ is equal to
$$
m_1=\sum_{k=1}^{n} {n\choose k} \frac{(-1)^{k+1}}{n!} \sum_{i=0}^{k- 1}
{n \choose i} (\alpha-1)^{n-i}(-\alpha)^i.
$$

Changing the order of summation and using the identity
\begin{equation}\label{th_2_1_2}
\sum_{k=0}^i (-1)^k {n\choose k}
= (-1)^i {{n-1}\choose i}
\end{equation}
(see, e.\,g., \cite{prudnikov_1981}) we get
 \begin{eqnarray}\label{th_2_1_3}
m_1 &=&
\frac{1}{n!} \sum_{i=0}^{n-1} {n \choose i} (\alpha-1)^{n-i}(-\alpha)^i
\sum_{k=0}^i (-1)^k {n\choose k} \nonumber\\
&=& \frac{1}{n!} \sum_{i=0}^{n-1} {n \choose i} { {n-1}\choose i}
  (\alpha-1)^{n-i}\alpha^i.
  \end{eqnarray}

Now, let us turn to $E^{(2)}.$ First, note that $E^{(2)}$
contains the domain
$S=\{x_i\geq 0,$ $ \sum x_i \leq 1\}$ with the measure  
$1/n!$. Next,
fix $k\in\{1,...,n\}$ and consider the subset
$G^\prime \subset E^{(2)}$
corresponding to the inequalities $x_1,\ldots,x_k<0;$ $x_{k+1},\ldots,x_n\geq
0.$ Put $y_1=-x_1,\ldots,y_k=-x_k;$ $y_{k+1}=x_{k+1},\ldots,y_n=x_n.$ Then
$$
G^\prime=\{y: y_{k+1}+\ldots+y_n\leq 1+ y_1+\ldots+y_k\leq \frac{\gamma-1}{2}, \
y_i\geq 0 \}.
$$

Let $\beta=({\gamma-1})/{2}.$ Then the following equalities hold:
\begin{eqnarray}
\mes(G^\prime) &=&
\int\limits_{0}^{\beta} \,dy_1
\int\limits_{0}^{\beta-y_1} \,dy_2
\ldots
\int\limits_{0}^{\beta-y_1-\ldots-y_{k-1}} \,dy_k
\nonumber\\ \nonumber\\
& &
\int\limits_{0}^{1+y_1+\ldots+y_k} \,dy_{k+1}
\int\limits_{0}^{1+y_1+\ldots+y_k-y_{k+1}} \,dy_{k+2}
\ldots
\int\limits_{0}^{1+y_1+\ldots+y_k-y_{k+1}-\ldots-y_{n-1}} \,dy_n \nonumber\\ \nonumber \\
&=&
\int\limits_{0}^{\beta} \,dy_1
\int\limits_{0}^{\beta-y_1} \,dy_2
\ldots
\int\limits_{0}^{\beta-y_1-\ldots-y_{k-1}} \frac{(1+y_1+\ldots+y_k)^{n-k}}{(n-k)!}
\,dy_k \nonumber\\ \nonumber\\
&=& \left[\sum_{j=0}^{k-1} (-1)^{k-1-j} \frac{(1+\beta)^{n-j}\beta^j}{(n-j) !j!}\right]
+\frac{(-1)^k}{n!}
= \frac{(-1)^{k+1}}{n!}
\left(\left[\sum_{j=0}^{k-1} {n\choose j} (1+\beta)^{n-j}(-\beta)^j\right] - 1 \right).\nonumber
\end{eqnarray}


Clearly, the set $E^{(2)}$ is the union of all such sets $G^\prime$
corresponding to various $k=1,\ldots,
n,$ and~the~simplex $S.$ Therefore,
$$
m_2 =\mes(E^{(2)})=
\frac{1}{n!} \left(\left\{ \sum_{k=1}^{n} (-1)^{k+1} {n\choose k}
\left(\left[\sum_{j=0}^{k-1} {n\choose j} (1+\beta)^{n-j}(-\beta)^j\right] - 1 \right)\right\}
+1 \right). 
$$

\par Note that
$1+\beta=(\gamma+1)/2=\alpha$ and $\beta=(\gamma-1)/2=\alpha-1$.
Let us make the substitution $i=n-j$ in the internal sum in the right hand side of this equality, and make use of the formula
$$\sum_{k=0}^{n} (-1)^k {n\choose k} = \sum_{k=1}^{n} (-1)^k
{n\choose k} + 1 = 0.
$$

\par We obtain the following:
\begin{eqnarray}
m_2 &=&
\frac{1}{n!} \Bigl(1+\sum_{k=1}^{n} (-1)^{k+1} {n\choose k}
\Bigl( (-1)^n\sum_{i=n-k+1}^{n} {n\choose i} (\alpha-1)^{n-i}(-\alpha)^i
  - 1\Bigr)
\Bigr) \nonumber\\
&=&
\frac{(-1)^n}{n!} \sum_{k=1}^{n} (-1)^{k+1} {n\choose k}
\sum_{i=n-k+1}^{n} {n\choose i} (\alpha-1)^{n-i}(-\alpha)^i.
\nonumber
\end{eqnarray}

Changing the order of summation, we obtain
$$
m_2 =
\frac{(-1)^n}{n!}
\sum_{i=1}^{n} {n\choose i} (\alpha-1)^{n-i}(-\alpha)^i
\sum_{k=n+1-i}^{n} (-1)^{k+1} {n\choose k}.$$
Using \eqref{th_2_1_2}, we can write
$$
\sum_{k=n+1-i}^{n} (-1)^{k+1} {n\choose k}
=\sum_{k=n+1-i}^{n} (-1)^{k+1}{n\choose {n-k}}
=\sum_{j=0}^{i-1} (-1)^{n-j+1}{n\choose j}=
(-1)^{n+i} { {n-1}\choose {i-1} }.$$
Therefore,
\begin{equation}\label{th_2_1_4}
m_2=
\frac{1}{n!}\sum_{i=1}^{n} {n\choose i} { {n-1}\choose {i-1} }
  (\alpha-1)^{n-i}
\alpha^i.
\end{equation}

\par From this, \eqref{th_2_1_3} and \eqref{th_2_1_4}, we have
\begin{eqnarray}
\mes(E_{n,\gamma})&=& m_1+m_2 \nonumber\\
&=& \frac{1}{n!} \sum_{i=0}^{n-1} {n \choose i} { {n-1}\choose i}
  (\alpha-1)^{n-i}\alpha^i +
\frac{1}{n!}\sum_{i=1}^{n} { n\choose i} { {n-1}\choose {i-1} }
  (\alpha-1)^{n-i}
\alpha^i \nonumber\\
&=& \frac{1}{n!} \sum_{i=1}^{n-1} { n \choose i}
\left( { {n-1}\choose i} + { {n-1}\choose {i-1} } \right)
  (\alpha-1)^{n-i}\alpha^i +
\frac{1}{n!} \Bigl((\alpha-1)^n+\alpha^n \Bigr)\nonumber\\ &=&
\frac{1}{n!}\sum_{i=0}^{n} {n\choose i}^2 (\alpha-1)^{n-i}
\alpha^i =
\frac{1}{2^n n!}\sum_{i=0}^{n} {n\choose i}^2 (\gamma-1)^{n-i}
(\gamma+1)^i \nonumber
\end{eqnarray}
completing the proof of the left-hand equality in \eqref{th_2_1_1}.
We took into account that $$\displaystyle { {n-1}\choose i } + { {n-1}\choose {i-1} } =
{n\choose i}.$$

The right-hand equality in \eqref{th_2_1_1} follows from the identity
$$
\sum_{i=0}^n {n\choose i}^2 t^i=(1-t)^n \chi_n\left(\frac{1+t}{1-t}\right).~~~~~
\text{See \cite{prudnikov_1981}.}
$$

\par Let us set $t={(\gamma - 1)}/{(\gamma+1)}$. Then
$(1-t)^n=2^n (\gamma+1)^{-n}$ and $(1+t)/(1-t)=\gamma$.
Hence,
\begin{eqnarray}
\mes(G)&=&
\frac{1}{2^n n!} \sum_{i=0}^{n} {n \choose i}^2 (\gamma - 1)^{n-i}
(\gamma + 1)^i=
\frac{1}{2^n n!} \sum_{i=0}^{n} {n \choose i}^2 (\gamma+1)^{n-i} \nonumber\\
&=& \frac{1}{2^n n!} (\gamma+1)^n
\sum_{i=0}^{n} {n \choose i}^2 \Bigl( \frac{\gamma-1}
{\gamma +1}\Bigr)^i =
\frac{\chi_n(\gamma)}{n!}. \nonumber
\end{eqnarray}

The  proof of Theorem \ref{th_mes_E_xi_eqs} is complete.
\hfill$\Box$

Let us give some simple examples. The set
$E_{1,2}=\{ x\in {\mathbb R}: |x|+|1-x|\leq 2\}$ is the segment 
$\left[-{1}/{2}, {3}/{2}\right]$ with the length ${\rm mes_1}(E_{1,2})=\chi_1(2)/1!=2$. The set
$E_{2,2}=\{ x\in {\mathbb R}^2: |x_1|+|x_2|+|1-x_1-x_2|\leq	 2\}$ is~a~hexagon on the plane with the area ${\rm mes_2}(E_{2,2}) =\chi_2(2)/2!=11/4$.



\par
\begin{figure}[h!]
\center{\includegraphics[scale=0.35]{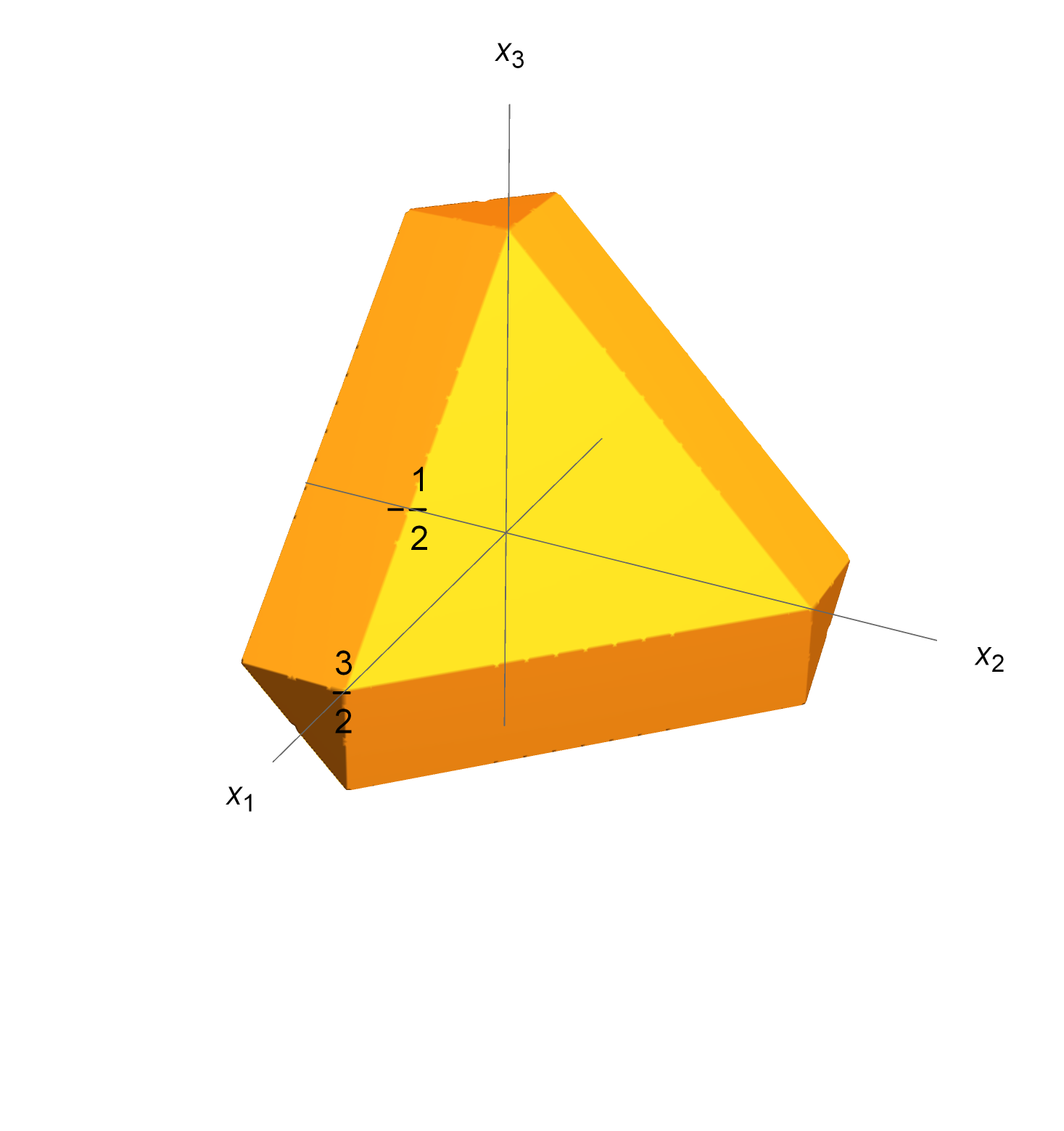}}
\vspace{-3cm}
\caption{The set $E_{3,2}$}
\label{fig:nev_ukl_example_E_3_2_one}
\end{figure}

 The three-dimensional domain
$$
E_{3,2}=\{ x\in {\mathbb R}^3: |x_1|+|x_2|+|x_3|+|1-x_1-x_2-x_3|\leq	 2\}
$$
is shown in Fig.~\ref{fig:nev_ukl_example_E_3_2_one}. It has the volume ${\rm mes_3}(E_{3,2})=\chi_3(2)/3!=17/6.$

\smsk


\par There is an interesting open problem  related to the above-mentioned equality \eqref{th_2_1_1}. Essentially,
this~equality
along with Rodrigues' formula and other well-known relations  gives a characterization of Legendre polynomials. This characterization is written via the volumes of some convex polyhedra. Namely, for $t\geq 1$
\begin{equation}\label{Legendre}
\chi_n(t)=
n!\,{\rm mes}_n(E_{n,t}),
\end{equation}
where $E_{n,t}$ is defined by (\ref{ENG}).

\par The question arises about analogues of \eqref{Legendre} for other classes of orthogonal polynomials, such as
Chebyshev polynomials or, more generally, Jacobi polynomials: {\it Is the equality \eqref{Legendre} a particular case of a more general pattern?}
We would be grateful for any information on this matter.

Let us note that, from  \eqref{Legendre} and \eqref{reccurent_Legendre}, we have
$$
{\rm mes}_{n+1}(E_{n+1,t})=\frac{2n+1}{(n+1)^2}\,t\, {\rm mes}_n(E_{n,t})-\frac{1}{(n+1)^2}{\rm mes}_{n-1}(E_{n-1,t}).
$$

The direct establishing this recurrence relation for the measures of $E_{n,t}$ could provide a proof
of~Theorem \ref{th_mes_E_xi_eqs} different from the above.


\SECT{5. Inequality  $\theta_n(K)
\geq
\chi_n^{-1}
\left(\dfrac{\vo(K)}{\simp(K)}\right)
$}{5}
\label{sect_theta_n_K_geq_xi_n_vol_K_simp_K}

\indent\par Based on Theorem \ref{th_mes_E_xi_eqs}, in this section we obtain lower bounds for the norm of the projection operator associated with linear interpolation on an arbitrary convex body $K$.

\begin{theorem}\label{th_norm_P_xi_K}
Suppose $P: C(K)\to \Pi_1\left({\mathbb R}^n\right)$ is an arbitrary interpolation
projector. Then for~the~corresponding simplex $S\subset K$ and the node matrix
${\bf A}$ we have
\begin{equation}\label{norm_P_xi_K_ineqs}
\|P\|_{K}
\geq
\chi_n^{-1} \left(\frac{n!\vo(K)}{|\det({\bf A})|}\right)=\chi_n^{-1} \left(\frac{\vo(K)}{\vo(S)}\right).
\end{equation}
\end{theorem}

\smallskip
{\it Proof.}
For each $i=1,$ $\ldots,$ $n$, let us  subtract from the
$i$th row of ${\bf A}$ its $(n+1)$th row.
Denote by ${\bf B}$ the submatrix of order $n$ which stands in the first
$n$ rows and columns of the result matrix. Then
$$
|\det({\bf B})|=|\det({\bf A})|= n!\vo(S)\leq n!\vo(K)
$$
so that
\begin{equation}\label{proof_111}
\frac{|\det({\bf B})|}{n!\vo(K)}\leq 1.
\end{equation}

Let $x^{(j)}$ be the vertices and $\lambda_j$ be the basic Lagrange polynomials of  simplex $S$.
Since $\lambda_1(x),$ $\ldots$ , $\lambda_{n+1}(x)$ are the barycentric coordinates of a point $x$ with respect to
$S$, we have
$$
\|P\|_{K}=\max_{x\in K} \sum_{j=1}^{n+1}|\lambda_j(x)|
=\max
\left\{ \sum_{j=1}^{n+1}|\beta_j|:
 \, \sum_{j=1}^{n+1}\beta_j=1,
\,
\sum_{j=1}^{n+1}\beta_jx^{(j)}\in K
\right\}.
$$

Let us replace $\beta_{n+1}$ with the equal value
$1-\sum\limits_{j=1}^n
\beta_j.$ The condition
$\sum\limits_{j=1}^{n+1}\beta_jx^{(j)}$ $\in K$ is equivalent to
$$
\sum\limits_{j=1}^{n}\beta_j(x^{(j)}-x^{(n+1)})\in K^\prime=K-x^{(n+1)}.
$$
Hence,
\begin{equation}\label{norm_P_barycentric_from1_ton}
\|P\|_{K}=\max\left \{ \sum_{j=1}^{n} \left|\beta_j\right| +
\left|1-\sum_{j=1}^n\beta_j\right|
\right \}
\end{equation}
where the maximum is taken over all
$\beta_j$ such that
$\sum\limits_{j=1}^n\beta_j(x^{(j)}-x^{(n+1)})
\in K^\prime.$ Clearly,
$
\vo(K^\prime)=\vo(K).
$

Let us consider the nondegenerate linear operator $F:{\mathbb R}^n \to {\mathbb R}^n,$
which maps the point $\beta=(\beta_1,\ldots,\beta_n)$ to  the point
$x=F(\beta)$
according to the rule $$x=\sum_{j=1}^n \beta_j \left(x^{(j)}-x^{(n+1)}\right).$$
We have the matrix equality
$F(\beta)=(\beta_1,\ldots,\beta_n) {\bf B},$
where $\bf B$ is the above introduced matrix of order $n$ with the entries $b_{ij}=x_j^{(i)}-x_j^{(n+1)}.$
Let us put
$$
\gamma^*=\chi_n^{-1}\left(\displaystyle \frac{n!\vo(K)}{|\det {\bf B}|}\right).
$$
Note that, thanks to \eqref{proof_111},
$(n!\vo(K))/|\det({\bf B})|\ge 1$ so that the constant $\gamma^*$ is well defined. Let us also note that $\chi_n(\gamma^*)=n!\vo(K)/|\det({\bf B})|$.
\smsk

\par Given $\gamma\geq 1$, let us introduce a set
$$
E_{n,\gamma}= \left\{ \beta=(\beta_1,\ldots,\beta_n)\in{\mathbb R}^n :
\sum_{j=1}^n |\beta_j| +
\left|1- \sum_{j=1}^n \beta_j\right|  \leq \gamma \right\}.
$$

Let us show that $K^\prime \not \subset F(E_{n,\gamma})$ provided $\gamma < \gamma^*.$ Indeed, this is immediate from Theorem \ref{th_mes_E_xi_eqs}  because, thanks to this theorem,
\begin{eqnarray}
\mes(F(E_{n,\gamma})) &<&
\mes(F(E_{n,\gamma^*}))=|\det {\bf B}|\cdot \mes(E_{n,\gamma^*}) \nonumber\\
&=&
|\det {\bf B}|\cdot
\frac{\chi_n(\gamma^*)}{n!}=\vo(K)=\vo(K^\prime). \nonumber
\end{eqnarray}

Thus, for every $\varepsilon > 0$ there exists a point
$x^{(\varepsilon)}$ with the properties:
$$x^{(\varepsilon)}=
\sum\beta_j^{(\varepsilon)} \left(x^{(j)}-x^{(n+1)}\right)\in K^\prime \mbox{\ \  and \  \ }
\left|\sum\beta_j^{(\varepsilon)}\right|
+\left|1-\sum\beta_j^{(\varepsilon)}\right|
\geq \gamma^*-\varepsilon.$$
In view of \eqref{norm_P_barycentric_from1_ton} this gives
$\|P\|_{K}\geq \gamma^*-\varepsilon.$
Since~$\varepsilon>0$ is an arbitrary, we obtain
$$\|P\|_{K} \geq \gamma^*=
\chi_n^{-1} \left(\frac{n!\vo(K)}{|\det({\bf B})|}\right)=
\chi_n^{-1} \left(\frac{n!\vo(K)}{|\det({\bf A})|}\right)=\chi_n^{-1} \left(\frac{\vo(K)}{\vo(S)}\right).$$
The theorem is proved.
\hfill$\Box$

\msk
Recall that $\simp(K)$ denotes the maximum volume of a nondegenerate simplex $S$ with vertices in~$K$
(see Definition \ref{simp_n_K_def}).

\begin{theorem}\label{th_theta_n_K_thru_vol_K_simp_K}
 Let $K$ be an arbitrary convex body in ${\mathbb R}^n$. Then
\begin{equation}\label{theta_n_K_ineq_thru_xi_minus_one}
\theta_n(K)
\geq
\chi_n^{-1}
\left(\frac{\vo(K)}{\simp(K)}\right).
\end{equation}
\end{theorem}

\smallskip
{\it Proof.} By \eqref{norm_P_xi_K_ineqs}, for any  interpolation
projector $P: C(K)\to \Pi_1\left({\mathbb R}^n\right)$,
$$
\|P\|_{K}
\geq
\chi_n^{-1} \left(\frac{\vo(K)}{\vo(S)}\right)\geq \chi_n^{-1}
\left(\frac{\vo(K)}{\simp(K)}\right), 
$$
where $S=S_P.$ This  immediately gives \eqref{theta_n_K_ineq_thru_xi_minus_one}. \hfill$\Box$

\smsk

\begin{remark}\label{RM-A}{\em If  simplex $S\subset K$ has  maximum volume, then $K\subset (n+2)S,$ see
\eqref{S_sub_K_sub_n_plus_2_S}. Therefore, the~ratio $\vo(K)/\simp(K)$ in
\eqref{theta_n_K_ineq_thru_xi_minus_one} is bounded from above by $(n+2)^n.$

In the case $K=Q_n$ we have $\vo(K)=1, \simp(K)=\nu_n$ (see Definition \ref{h_n_nu_n_def}). If $K=B_n$,  then $\vo(K)=\varkappa_n,
\simp(K)=\sigma_n$ (see Definition \ref{kappa_n_sigma_n_def}). Thus, in these cases, \eqref{theta_n_K_ineq_thru_xi_minus_one} imply the following two inequalities:
\begin{equation}\label{corol_2_1}
\theta_n(Q_n)\geq
\chi_n^{-1} \left(\displaystyle\frac{1}{\nu_n}\right)
~~~~~\text{and}~~~~~~
\theta_n(B_n)
\geq
\chi_n^{-1}
\left(\frac{\varkappa_n}{\sigma_n}\right).
\end{equation}
\par In the next section we will give two important  consequences of these inequalities.}\rbx
\end{remark}

\begin{remark}\label{RM-B}{\em
\par Inequalities  \eqref{xi_n_K_leq_n+2_theta_n_K_leq_n+1} and \eqref{theta_n_K_ineq_thru_xi_minus_one} show that for each convex body $K\subset \RN$ 
$$ \chi_n^{-1}
\left(\frac{\vo(K)}{\simp(K)}\right)\leq \theta_n(K)\leq n+1.
$$
 In Section 6 we prove analogues of these inequalities for an arbitrary (not necessarily convex) compact set $E\subset \RN$ with $\vo(\conv(E))>0$.}\rbx
\end{remark}


\SECT{6. Inequalities $\theta_n(Q_n)>c\sqrt{n}$ \  and
\ $\theta_n(B_n)>c\sqrt{n}$}{6}
\label{sect_theta_n_cube_ball_geq_sqrt_n}

\smallskip
The Stirling formula (see, e.\,g., \cite{fiht_2001})
$n!=\sqrt{2\pi n}\,(n/e)^n e^{\frac{\zeta_n}{12n}}$, $0<\zeta_n<1$, yields
\begin{equation}\label{n_fact_ineqs}
\sqrt{2\pi n}\left(\frac{n}{e}\right)^n<n!<\sqrt{2\pi n}\left(\frac{n}{e}\right)^n
e^{\frac{1}{12n}}.
\end{equation}

Let us also see that for every $n>1$ the following inequality
\begin{equation}\label{chi_ineq}
\chi_n^{-1}(s)>\left( \frac{s}{ {n\choose{\lfloor n/2\rfloor}} }\right)^{1/n}
\end{equation}
holds. In fact, if $t\geq 1$ and $n>1$, then, according to \eqref{th_2_1_1},
\begin{eqnarray}
\chi_n(t) &=&
\frac{1}{2^n}\sum_{i=0}^{n} {n\choose i}^2 (t-1)^{n-i}(t+1)^i
< {n\choose{\lfloor n/2 \rfloor}} \cdot \frac{1}{2^n}\sum_{i=0}^{n}
{n\choose i} (t-1)^{n-i}(t+1)^i \nonumber\\ \nonumber\\
&=& {n\choose{\lfloor n/2 \rfloor}}\cdot (2t)^n \cdot 2^{-n}=
{n\choose{\lfloor n/2 \rfloor} }\, t^n \nonumber
\end{eqnarray}
proving \eqref{chi_ineq}.
For even $n$, we have
$$
\displaystyle{n\choose{\lfloor n/2 \rfloor} }=
{n\choose{n/2}}=\frac{n!}{\left((n/2)!\right)^2},
$$
therefore,
\begin{equation}\label{chi_ineq_2}
\chi_n^{-1}(s)>\left( \frac{s\left((n/2)!\right)^2}
{ n!}\right)^{1/n}.
\end{equation}
\par If $n$ is odd, then
$${n\choose{\lfloor n/2 \rfloor}}
=\frac{n!}{  \frac{n+1}{2} ! \frac{n-1}{2}!  }$$
so that, thanks to \eqref{chi_ineq},

\begin{equation}\label{chi_ineq_3}
\chi_n^{-1}(s)>\left( \frac{s\,\frac{n+1}{2}!\frac{n-1}{2}!}
{ n!}\right)^{1/n}.
\end{equation}

\begin{theorem}\label{th_est_theta_n_cube_sqrt_n}
 For all $n\in {\mathbb N}$,
\begin{equation}\label{th_2_4_form}
\theta_n(Q_n)>\frac{\sqrt{n-1}}{e}.
\end{equation}
\end{theorem}

\smallskip
{\it  Proof.} The case $n=1$ is trivial since $\theta_1(Q_1)=1$.  If $n>1$, from
\eqref{chi_ineq_2}, \eqref{chi_ineq_3}
and the~Hadamard inequality \eqref{adamar_clements_lindstrem}, we have
$\nu_n\le \left(n+1\right)^{(n+1)/2}/2^nn!$.

\par For even $n$, thanks to the first inequality in \eqref{corol_2_1},
the left-hand inequality in \eqref{n_fact_ineqs}, and
  \eqref{chi_ineq_2},
we get
\begin{eqnarray}
\theta_n(Q_n)
&\geq&
\chi_n^{-1}\left(\frac{1}{\nu_n}\right)\geq
\chi_n^{-1}\left( \frac{2^nn!}{(n+1)^{(n+1)/2}}\right) >
2\left(\frac{[(n/2)!]^2}{(n+1)^{(n+1)/2}}\right)^{1/n}
\nonumber\\
&>&
\frac{2}{(n+1)^{1/2+1/(2n)}}\left( \sqrt{\pi n} \left(\frac{n}{2e}\right)^{n/2}
\right)^{2/n}
=
\frac{ \left(\pi n\right)^{1/n}n }{e(n+1)^{1/2+1/(2n)}}>
\frac{\sqrt{n-1}}{e}. \nonumber
\end{eqnarray}
If $n$ is odd, then, thanks to \eqref{chi_ineq_3},
\begin{eqnarray}
\theta_n(Q_n)
&\geq& \chi_n^{-1}\left(\frac{1}{\nu_n}\right)\geq
\chi_n^{-1}\left( \frac{2^nn!}{(n+1)^{(n+1)/2}}\right)
>
\left( \frac{2^n \frac{n+1}{2}!\frac{n-1}{2}! } {(n+1)^{(n+1)/2}}
\right)^{1/n}  \nonumber \\   \nonumber \\
&>&
2\left( \frac{ \pi\sqrt{n^2-1}\left(n^2-1\right)^{(n-1)/2}
(n+1) }{(2e)^n} \right)^{1/n}  \nonumber \\  \nonumber \\
&=&
\frac{1}{e}\pi^{1/n}\sqrt{n-1}(n+1)^{1/(2n)}>
\frac{\sqrt{n-1}}{e} \nonumber
\end{eqnarray}
proving that \eqref{th_2_4_form} holds for all $n$.
\hfill$\Box$

\smallskip
In some situations, the estimates of Theorem 2.4 can be improved \cite{nevskii_monograph}. In particular, we prove
that $\theta_n>\sqrt{n}/{e}$ if $n$ is even, and
$\displaystyle \theta_n>n/(e\sqrt{n-1})$ provided $n>1$ and $n\equiv 1({\rm mod}~4)$.

Let us also note that, thanks to \eqref{th_2_4_form}, the inequality $\theta_n(Q_n)> c\sqrt{n}$ holds with some $c$, $0<c<1$.
A~suitable estimate is
\begin{equation}\label{theta_n_lt_sqrt_n}
\theta_n(Q_n)>\frac{2\sqrt{2}}{3e} \, \sqrt{n}.
\end{equation}

Indeed, if $n\leq 8$, then the right-hand side of \eqref{theta_n_lt_sqrt_n} is less than 1,  while for
$n\geq 9$
$$ \theta_n(Q_n)>\frac{\sqrt{n-1}}{e}\geq \frac{2\sqrt{2}}{3e} \, \sqrt{n}.$$
Notice that
$\displaystyle\frac{2\sqrt{2}}{3e} = 0.3468 \ldots$
\bsk

This approach was extended to linear interpolation on the unit ball $B_n$, see~\cite{nevskii_mais_2019_26_3}.

\begin{corollary}\label{corol_est_theta_n_ball_thru_xi_n}
For every $n$,
\begin{equation}\label{theta_n_chi_n_through_n_ineq}
\theta_n(B_n)
\geq
\chi_n^{-1}
\left(\frac{\pi^{\frac{n}{2}}n!}{\Gamma\left({n}/{2}+1\right)\sqrt{n+1}
\left({(n+1)}/{n}\right)^{{n}/{2}}}
\right).
\end{equation}
If $n=2k$, then  (\ref{theta_n_chi_n_through_n_ineq}) is equivalent to the inequality
\begin{equation}\label{theta_n_chi_n_through_n_is_2k_ineq}
\theta_{2k}(B_{2k})
\geq
\chi_{2k}^{-1}
\left(\frac{\pi^{k}(2k)!}{k!\sqrt{2k+1}
\left({(2k+1)}/{2k}\right)^k}
\right).
\end{equation}
For $n=2k+1$ we have
\begin{equation}\label{theta_n_chi_n_through_n_is_2k_plus_1_ineq}
\theta_{2k+1}(B_{2k+1})
\geq
\chi_{2k+1}^{-1}
\left(\frac{2(k!)(4\pi)^{k}}{\sqrt{2k+2}
\left({(2k+2)}/{(2k+1)}\right)^{{(2k+1)}/{2}}}
\right).
\end{equation}
\end{corollary}

\smallskip
{\it Proof.} It is sufficient to apply
the second inequality in (\ref{corol_2_1}),
(\ref{kappa_n_sigma_n_whole_1}), and
(\ref{kappa_n_even_and_odd_1}).
\hfill$\Box$
\msk

Thanks to \eqref{chi_ineq_2} and \eqref{chi_ineq_3}, inequalities \eqref{theta_n_chi_n_through_n_is_2k_ineq} and
\eqref{theta_n_chi_n_through_n_is_2k_plus_1_ineq} imply the following low bound for the constant $\theta_n(B_n)$
(see  \cite{nevskii_mais_2021_28_2}).

\begin{theorem}\label{th_theta_n_ball_geq_c_sqrt_n}
 There exists an absolute constant $c>0$ such that
\begin{equation}\label{theta_n_B_n_gt_c_sqrt_n_1}
\theta_n(B_n)>c\sqrt{n}.
\end{equation}

Inequality
(\ref{theta_n_B_n_gt_c_sqrt_n_1}) takes place  with the constant

\begin{equation}\label{C1}
c=
\frac{  \sqrt[3]{\pi} }{\sqrt{12e}\cdot\sqrt[6]{3} } =  0.2135...\,.
\end{equation}
\end{theorem}
\bsk


\SECT{7. Estimates for linear interpolation on an arbitrary compact set}{7}\label{sect_Est_on_a_compact}

\indent
\par In this section we show how the above approach leads to estimates of the minimum projector norm under linear interpolation on a  compact set.

Let $E$ be an arbitrary compact in $\RN$. Everywhere in this section, we let $K$ denote {\it the convex hull of the set $E$}. We will assume that $\vol(K)>0$. By $\|P\|_E$ we denote the norm of the interpolation projector
$P:C(E)\to \Pi_1(\RN)$ as an operator from $C(E)$ to $C(E)$. We let $\theta_n(E)$
denote the minimal norm  $\|P\|_E$ over all projectors $P$ with nodes in $E$. By $\simp(E)$ we denote
the maximum volume of a simplex with vertices in $E.$

\smsk
\par We will need the following elementary lemma.

\begin{lemma}\label{lemma_conv_func_on_Rn}
If $\varphi:\RN\to {\mathbb R}$ is a convex continuous function, then $\max\limits_K \varphi=\max\limits_E \varphi.$
\end{lemma}

\smsk
{\it Proof.} The maxima from the condition exist since $\varphi$ is a continuous function. As $K=\conv(E)$, for~any
$y\in K$ there exist a number $m$,  points $y^{(1)},
\ldots, y^{(m)}\in E$ and numbers $\mu_1,\ldots,\mu_m$ such that
$$y=\sum_{i=1}^m \mu_i y^{(i)}, \quad \mu_i\geq 0, \ \mbox{ \ and \ } \sum_{i=1}^m \mu_i=1.$$
Clearly, $\varphi(y^{(i)})\leq  \max\limits_E \varphi.$ The convexity of $\varphi$ implies 
$$\varphi(y)=\varphi \left(\sum_{i=1}^m \mu_i y^{(i)}\right) \leq \sum_{i=1}^m \mu_i \varphi(y^{(i)})\leq
\left(\sum_{i=1}^m \mu_i\right) \max\limits_E \varphi = \max\limits_E \varphi.$$
Therefore, $ \max\limits_K\varphi\leq\max\limits_E \varphi$. The inverse  inequality is obvious. \hfill$\Box$

\smsk
Clearly, the result of Lemma \ref{lemma_conv_func_on_Rn} is well-known.  It is immediate from
the following Bauer's maximum principle \cite{bauer_1958}.
{\it Any function that is convex and continuous, and defined on a convex and compact set, attains its maximum at some extreme point of that set.} Consequently, 
maximum $\varphi$ on $K=\conv(E)$  is attained at some extreme point of $E$.

\smsk
\begin{theorem}\label{theor_theta_n_E_thru_vol_K_simp_E}
 We have
\begin{equation}\label{theta_n_E_ineq_thru_xi_minus_one}
\theta_n(E)
\geq
\chi_n^{-1}
\left(\frac{\vo(K)}{\simp(E)}\right).
\end{equation}
\end{theorem}

\smsk
{\it Proof.} First, let us note that for an arbitrary polynomial $p\in \Pi_1(\RN)$ the following inequality
\begin{equation}\label{norm_p_E_eq_norm_p_p_K}
\|p\|_E=\|p\|_K
\end{equation}
holds. (Recall that $K=\conv(E)$.) This is immediate from Lemma \ref{lemma_conv_func_on_Rn} for the convex  continuous function
$\varphi(x)=|p(x)|.$

\par Let $P:C(E)\to \Pi_1(\RN)$ be an arbitrary interpolation projector with nodes in $E$. We will consider it also as a projector from $C(K)$. Thanks to \eqref{norm_p_E_eq_norm_p_p_K}, it follows that
$\|P\|_E=\|P\|_K.$
Therefore, thanks to inequality \eqref{theta_n_K_ineq_thru_xi_minus_one}
of Theorem \ref{th_theta_n_K_thru_vol_K_simp_K},
$$
\theta_n(E)\geq \theta_n(K)
\geq \chi_n^{-1}
\left(\frac{\vo(K)}{\simp(K)}\right).
$$

\par Finally, let us see that $\simp(K)=\simp(E)$. Indeed,
let us suppose $S\subset K$ is any simplex with some vertex $x\not\in E$. Without changing other vertices, we can replace vertex $x$  by vertex $x^\prime$ so that the volume of the resulting simplex does not decrease. Indeed, this volume increases with $\dist(x; \Gamma)$, where $\Gamma$ is the $(n-1)$-dimensional hyperplane containing all the vertices of the simplex except $x$. Let~$\Gamma$ be given by the equation $q(z)=\langle a, z\rangle + a_0=0$,
$a=(a_1,\ldots,a_n)\in \RN, a_0\in {\mathbb R}$. Then
$$\dist(x; \Gamma)=\frac{|q(x)|}{||a|||}$$
is  obviously a convex continuous function. By Lemma \ref{lemma_conv_func_on_Rn}, the maximum of
$\dist(x; \Gamma)$ 
on $K$ is attained at some point $x^\prime\in E.$ Applying this procedure sequentially to all vertices of the simplex not belonging to $E$, we construct a new simplex with vertices on $E$ without reducing the volume.

\par Thus, $\simp(K)=\simp(E)$, and the proof  is complete.\hfill$\Box$

\smsk
\begin{theorem}\label{theor_theta_n_E_leq_nplus1}
Let $E\subset\RN$ be a compact set such that $\vo(\conv(E))>0$. Then the following inequality
\begin{equation}\label{theta_n_E_leq_nplus1}
\theta_n(E) \leq n+1
\end{equation}
holds.
\end{theorem}

\smsk
{\it Proof.} For an arbitrary interpolation projector $P:C(E)\to \Pi_1(\RN)$, we have the following equality
$$
\|P\|_E=\max_{x\in E} \sum_{j=1}^{n+1} |\lambda_j(x)|.
$$
This can be proved in the same way as in the case  
$E=K$ (see \eqref{norm_P_intro_cepochka}).

\par Let $S$ be a simplex with maximum volume over all family of simplices with vertices in $E$. For this simplex,
 $|\lambda_j(x)|\leq 1$
(see the proof of
Theorem \ref{th_norm_P_for_max_simplex}, inequality (\ref{LM-1})). By the above formula, the corresponding projector $P=P_S$
satisfies $\|P_S\|_E\leq n+1$. Consequently, $\theta_n(E) \leq n+1$.
\hfill$\Box$

\smsk
Let us combine inequalities  \eqref{theta_n_E_ineq_thru_xi_minus_one} and \eqref{theta_n_E_leq_nplus1}.
 If $E\subset  \RN$ is an arbitrary compact set satisfying the~condition $\vo(\conv(E))>0$, then
$$
\chi_n^{-1}
\left(\frac{\vo(\conv(E))}{\simp(E)}\right)\leq \theta_n(E) \leq n+1.$$


\SECT{8. Concluding remarks and open questions
}{8}
\label{sect_Concluding_remarks}

\indent\par Despite the apparent simplicity of formulation,
 the problem of finding exact values of $\theta_n(Q_n)$ is~very difficult.
Since 2006, these values are known only for $n=1,2,3,$ and $7$ (see \cite{nevskii_mais_2006_13_2},
\cite{nevskii_monograph}) . Namely,
$$\theta_1(Q_1)=1, \quad \theta_2(Q_2)=\frac{2\sqrt{5}}{5}+1=1.8944\ldots, \quad \theta_3(Q_3)=2, \quad \theta_7(Q_7)= \frac{5}{2}.$$
The corresponding numbers $\xi(Q_n)$ are
$$
\xi_1(Q_1)=1, \quad \xi_2(Q_2)=\frac{3\sqrt{5}}{5}+1=2.3416\ldots, \quad \xi_3(Q_3)=3, \quad \xi_7(Q_7)=7.
$$
Hence, for each $n=1,2,3,7$ the right-hand inequality in  \eqref{nev_ksi_n_K_theta_n_K_ineq}
becomes an equality:
\begin{equation}\label{righthand_eq_theta_ksi_cube}
\xi_n(Q_n) =
\frac{n+1}{2}\left( \theta_n(Q_n)-1\right)+1.
\end{equation}

If a  nondegenerate simplex $S$ is contained in $Q_n$, then $d_i(S)\leq 1$. Using \eqref{alpha_d_i_S_eq}
we get
$$
\xi(Q_n;S)\geq \alpha(Q_n;S)=\sum_{i=1}^n \frac{1}{d_i(S)}\geq n.
$$
Therefore, always $\xi_n(Q_n)\geq n$.
If $n+1$  is an Hadamard number, then $\xi_n(Q_n)= n$ (see \cite{nevskii_dcg_2011,nevskii_monograph}). 

Thanks
to  \eqref{nev_ksi_n_K_theta_n_K_ineq}, inequality $\xi_n(Q_n)\geq n$ implies
\begin{equation}\label{theta_n_cube_geq_3_minus_frac}
\theta_n(Q_n)\geq 3-\frac{4}{n+1}.
\end{equation}
If  $n=1, 3,$ or $7$, then in \eqref{theta_n_cube_geq_3_minus_frac} we have an equality.

For $1\leq n\leq 3$  the  simplices corresponding to minimal projectors are just the same that the extremal simplices with respect to $\xi_n(Q_n)$. The proofs are given in \cite{nevskii_mais_2009_16_1,nevskii_monograph} (see also the recent survey \cite{nevskii_arxiv_2024}). Let us mention here these results.

The case $n=1$ is trivial: $\theta_1(Q_1)=1$ and a unique extremal simplex is the segment $[0,1]$ coinciding with $Q_1$.
Suppose $n=2.$
Up to rotations, the only  simplex extremal with respect both to $\theta_2(Q_2)$
and $\xi_2(Q_2)$ is the triangle with vertices $(0,0),$ $(1,\tau),$
$(\tau,1)$, where
$\tau=(3-\sqrt{5})/2=0.3819\ldots$ .
 This number
satisfies
$\tau^2-3\tau+1=0$ or
$$\frac{\tau}{1-\tau}=  \frac{1-\tau}{1}.$$
Hence,  $\tau$ delivers {\it the golden section} of the segment [0,1].
Sharp inequality
$\|P\|_{Q_2}\geq 2\sqrt{5}/5+1$
for projectors corresponding to simplices $S\subset Q_2$ gives
the new characterization of this classical notion.

Up to coordinate substitution,   each extremal simplex in $Q_3$ coincides with
the tetrahedron $S^\prime$ with vertices
$$(1,1,0),~~~(1,0,1),~~~(0,1,1),~~~(0,0,0)$$ or
the tetrahedron $S^{\prime\prime}$ with vertices
$$
\left(\frac{1}{2},0,0\right),~~~
\left(\frac{1}{2},1,0\right),~~~
\left(0,\frac{1}{2},1\right),~~~
\left(1,\frac{1}{2},1\right).
$$
In other words, if $\|P_S\|_{Q_3}=3,$ then  either vertices of $S$ coincide with vertices of the cube and~form
a regular tetrahedron or coincide with the centers of opposite edges of two opposite faces of the cube and
does not belong to a common plane.

Let us turn to the case $n=7$. Since $8$ is an Hadamard number, there exists a seven-dimensional regular simplex having vertices at vertices of the cube.
We can take the simplex with the vertices
$$
(1,1,1,1,1,1,1),~~~(0,1,0,1,0,1,0),~~~(0,0,1,1,0,0,1),~~~
(1,0,0,1,1,0,0),$$
$$
(0,0,0,0,1,1,1),~~~(1,0,1,0,0,1,0),~~~(1,1,0,0,0,0,1),~~~
(0,1,1,0,1,0,0).
$$

The equality $\xi_7(Q_7) =7$  and
inequality
$$
\xi_n(Q_n) \leq \frac{n+1}{2}\left( \theta_n(Q_n)-1\right)+1
$$
imply $\theta_7(Q_7) \geq  {5}/{2}$. But for the corresponding projector, $\|P\|_{Q_7}={5}/{2}$.
Therefore, $\theta_7 = {5}/{2}$, and this projector is minimal.

 Let $n+1$ be an Hadamard number, and let $S$ be an $n$-dimensional regular simplex having the vertices at vertices of  $Q_n$
 Then, for the
corresponding  projector
$P_S:C(Q_n)\to \Pi_1({\mathbb R}^n)$, we have
\begin{equation}\label{th2_formula}
\|P_S\|_{Q_n}\leq \sqrt{n+1}.
\end{equation}
Various proofs of \eqref{th2_formula} are given in \cite{nevskii_mais_2003_10_1}  and \cite{nevskii_monograph}.
Paper  \cite{nevskii_mais_22} contains the proof essentially making use of the structure of an Hadamard matrix.
The equality
 $\|P_S\|_
 {Q_n}
 =\sqrt{n+1}$ may hold as for all regular simplices $S$ having vertices at vertices of~the~cube $(n=1, n=3)$, as for some of them $(n=15),$
 or  may not be executed at all.

As it is shown in Section 5,
for each $n$,
$\theta_n(Q_n)\geq \sqrt{n-1}/e.$ Therefore,
if $n+1$ is an Hadamard number, then
$$\frac{\sqrt{n-1}}{e}\leq \theta_n(Q_n)\leq  \sqrt{n+1}.$$
In other words, the above estimate $\theta_n(Q_n)\geq c\sqrt{n}$ is sharp at least when $n+1$ is an Hadamard
number. For these dimensions, $\theta_n(Q_n)\asymp \sqrt{n}$.

\msk
The upper bounds of $\theta_n(Q_n)$ for special $n$ were improved by
A. Ukhalov and his students
\linebreak with the~help of
computer methods. In particular,  simplices of maximum volume in the cube were considered.
In all situations where $n+1$ is an Hadamard number, the full set of Hadamard matrices of the corresponding order was used. In particular, to obtain an estimate for $\theta_{23}$, all existing $60$ Hadamard matrices of order $24$ were considered. To estimate $\theta_{27}$,
we have to consider $487$ Hadamard matrices of order $28$.
Known nowaday upper estimates for $1\leq n\leq 27$ are given in \cite{nev_ukh_posobie_2022}. Here they are (for~brevity, we write $\theta_n=\theta_n(Q_n)$):
$$\theta_1=1, \quad \theta_2=\frac{2\sqrt{5}}{5}+1, \quad  \theta_3=2,
\quad \theta_4\leq \frac{3  (4+\sqrt{2})}{7},\quad\theta_5\leq  2.448804,$$
$$\theta_6\leq 2.6000\ldots ,\quad \theta_7=\frac{5}{2}, \quad \theta_8\leq \frac{22}{7},\quad\theta_9\leq 3, \quad\theta_{10}\leq \frac{19}{5}, \quad\theta_{11}\leq 3,$$
$$\theta_{12}\leq  \frac{17}{5},\quad
\theta_{13}\leq  \frac{49}{13},\quad \theta_{14}\leq \frac{21}{5},
\quad \theta_{15}\leq \frac{7}{2}, \quad \theta_{16}\leq  \frac{21}{5},\quad \theta_{17}\leq \frac{139}{34},
$$
 $$ \theta_{18}\leq 5.1400\ldots,
\quad \theta_{19}\leq 4, \quad \theta_{20}\leq 4.68879\ldots, \quad \theta_{21}\leq
\frac{251}{50},
\quad\theta_{22}\leq  \frac{1817}{335},$$
$$\theta_{23}\leq  \frac{9}{2},
\quad \theta_{24}\leq  \frac{103}{21}, \quad \theta_{25}\leq  5,
\quad \theta_{26}\leq  \frac{474}{91},
\quad \theta_{27}\leq  5.$$


The best  nowaday known  lower bound of $\theta_n(Q_n)$ has the form
 \begin{equation}\label{nev_theta_n_max_ineq}
\theta_n(Q_n)
\geq  \max \left[
3 - \frac{4}{n+1}, \,
\chi_n^{-1} \left(\frac{1}{\nu_n}\right) \right].
\end{equation}
Here $\chi_n$ is the standardized Legendre polynomial of degree $n$, see Section 3.
The values of the right-hand side of \eqref{nev_theta_n_max_ineq} for $1\leq n\leq 54$ are given in \cite{nev_ukh_posobie_2022}.

As noted in Section 1 (see  \eqref{nev_ksi_P_ineq}),
the inequality
\begin{equation}\label{nev_ksi_n_teta_n_ineq_sec4}
\xi_n(Q_n) \leq
\frac{n+1}{2}\left( \theta_n(Q_n)-1\right)+1
\end{equation}
 is true for any $n$.
So far, we  know only four values of $n$ for which this relation becomes an
equality: $n=1, 2, 3$, and $7$. These are exactly the cases in which
we know the exact values both of $\theta_n(Q_n)$ and $\xi_n(Q_n)$. In
\cite{nev_ukh_mais_2016_23_5}
the authors conjectured that the minimum of $n$ for which inequality (\ref{nev_ksi_n_teta_n_ineq_sec4}) is strict is $4.$ This is still an open problem.

The  above given estimate $\xi_n(Q_n)\geq n$ occurs to be exact in order of $n$.
 If $n>2$, then
\begin{equation}\label{xi_n_leq_frac_general}
\xi_n(Q_n)\leq \frac{n^2-3}{n-1}
\end{equation}
(see  \cite{nevskii_mais_2011_2,nevskii_monograph}).
If $n>1$, the right-hand side of  \eqref{xi_n_leq_frac_general} is strictly smaller than  $n+1$.
Inequality $\xi_n<n+1$
holds true also for $n=1,2$.
Thus, always $n\leq \xi_n(Q_n)< n+1$, i.e., $\xi_n(Q_n)-n\in [0,1)$.
However,  the exact values of the constant
$\xi_n(Q_n)$ are currently only known for
$n=2$, $n=5$ and $n=9$, as well as for an infinite set of $n$ for which there exists an Hadamard matrix of order $n+1$.

In all these cases, except $n=2$, the equality $\xi_n(Q_n)=n$ holds.
In the noted Hadamard's case, one can give the proof using the structure of Hadamard matrix of order $n+1$, see \cite{nev_ukh_beitrage_2018}.
In \cite{nev_ukh_beitrage_2018}, we have also discovered the exact values of $\xi_n(Q_n)$ for $n=5$ and $n=9$ and  constructed several infinite families of extremal simplices for $n=5,7,9$.

\par Thanks to the equivalence $\xi_n(Q_n)\asymp n$ and inequality $\theta_n\geq c \sqrt{n}$, for all sufficiently large $n$, we have
\begin{equation}\label{nev_strict}
\xi_n(Q_n) <
\frac{n+1}{2}\left(\theta_n(Q_n)-1\right)+1.
\end{equation}

Let $n_0$ be the minimal natural number
such that for all $n\geq n_0$ inequality (\ref{nev_strict}) holds. The problem about the exact value of $n_0$ is very difficult. The known
lower and upper bounds differ quite significantly.
From the preceding, we have the estimate $n_0\geq 8.$
In 2009 we proved that $n_0\leq 57$ (see \cite{nevskii_mais_2009_16_1,nevskii_monograph}).
A sufficient condition for the validity of
(\ref{nev_strict})
for $n>2$ is the inequality
\begin{equation}\label{nev_chi_nu_ineq}
\chi_n \left(\frac{3n-5}{n-1}\right)\cdot\nu_n<1.
\end{equation}

It was proved  in \cite{nevskii_mais_2009_16_1} that \eqref{nev_chi_nu_ineq}  is satisfied for $n\leq 57$.
Later calculations allowed the upper bound of $n_0$ to be slightly lowered.
 In \cite{nev_ukh_mais_2018_25_3} 
   it is noted
that $n_0\leq 53$. In~other words,
inequality \eqref{nev_strict} is satisfied at least
starting from $n=53$. A better estimate from above for $n_0$ is an open problem.
\smsk

\par Let us now proceed to interpolation on the unit ball $B_n$.

\par It was shown in \cite{nevskii_mais_2018_25_6} that  $\xi_n(B_n)=n$ for all $n$. Moreover, for a simplex
$S\subset B_n$, the equality $\xi(B_n;S)=n$ is equivalent to the fact that $S$ is a regular simplex inscribed in the ball.
By \eqref{nev_ksi_P_ineq}, for any projector $P:C(B_n)\to\Pi_1\left({\mathbb R}^n\right)$
we have
 \begin{equation}\label{nev_ksi_P_ball_ineq}
 \|P\|_{B_n}\geq 3- \frac{4}{n+1}.
 \end{equation}
 The right-hand equality in \eqref{nev_ksi_n_K_theta_n_K_ineq}, i.\,e., the equality
 \begin{equation}\label{nev_ksi_n_K_theta_n_ball_eq}
\xi_n(B_n) =
\frac{n+1}{2}\Bigl( \theta_n(B_n)-1\Bigr)+1,
\end{equation}
is equivalent to
\begin{equation}\label{theta_n_ball_eq_3_minus_frac}
 \theta_n(B_n)= 3- \frac{4}{n+1}.
 \end{equation}
As it shown in \cite{nev_ukh_mais_2019_26_2}, equalities
\eqref{nev_ksi_n_K_theta_n_ball_eq} -- \eqref{theta_n_ball_eq_3_minus_frac} take place for
$1\leq n\leq 4,$ while starting from $n=5$ we have the strict inequality
$$\xi_n(B_n) <
\frac{n+1}{2}\Bigl( \theta_n(B_n)-1\Bigr)+1.$$
For $1\leq n\leq 4,$ equality in \eqref{nev_ksi_P_ball_ineq} holds only if $S_P$ is a regular inscribed simplex.
\smsk

\par The fact that the lower estimate $\theta_n(B_n)>c\sqrt{n}$ obtained in Section 5 is exact in dimension $n$ was first established in  \cite{nev_ukh_mais_2019_26_2}: {\it the equivalence $\theta_n(B_n)\asymp \sqrt{n}$ is valid.}

\par A complete solution to the problem about the values of the numbers $\theta_n(B_n)$ was given in
\cite{nevskii_matzam_23}. Let us briefly describe these results. Let $\psi:[0,n+1]\to\R$ be a function defined by
$$\psi(t)=\frac{2\sqrt{n}}{n+1}\Bigl(t(n+1-t)\Bigr)^{1/2}+
\left|1-\frac{2t}{n+1}\right|,$$
and let
$$
a=a_n=\left\lfloor\frac{n+1}{2}
-\frac{\sqrt{n+1}}{2}\right\rfloor.
$$
\par Let $S$ be a regular simplex inscribed in $B_n$, and let $p_n$ be the norm of $P_S$. It was proved in
 \cite{nev_ukh_mais_2019_26_2} that

$$p_n=\max\{\psi(a),\psi(a+1)\} \mbox{\  \ and  \ }  \sqrt{n}\leq p_n \leq \sqrt{n+1}.$$

\smsk
\noindent Moreover, $p_n=\sqrt{n}$ \, only for $n=1$ and $p_n=\sqrt{n+1}$ holds if and only if
$\sqrt{n+1}$ is an integer.
\smsk

The equality $\theta_n(B_n)=p_n$ was obtained first for $1\leq n\leq 4$ (different proofs
are given in \cite{nev_ukh_mais_2019_26_2} and
\cite{nevskii_mais_2021_28_2}). As a conjecture for all $n$, this statement
was formulated in \cite{nevskii_mais_2021_28_2}. In \cite{nevskii_matzam_23} we developed a certain new  geometric approach to the problem which enabled us to prove this conjecture, i.e., the equality
$\theta_n(B_n)=p_n$ for an arbitrary positive integer $n$.

\par Thus, the following ineuality
$\sqrt{n}\leq \theta_n(B_n) \leq \sqrt{n+1}$ holds for all $n$. Furthermore, we proved that the minimal is any projector corresponding to a regular simplex inscribed into the boundary sphere, and there are no other minimal projections in this case.

 Let $k_n$ coincide with that of the numbers $a_n$ and $a_n+1$ on which $\psi(t)$ takes a larger value.
 The numbers $k_n$  increase with $n$, but not strictly monotonically. If $n\geq 2$, then
$k_n\leq n/2.$
As an example, we give below the numbers $k_n$ for $1\leq n\leq 15$, $n=50, n=100$, and $n=1000$ (\cite{nev_ukh_mais_2019_26_2}):
$$k_1=k_2=k_3=k_4=1, \quad k_5=k_6=2, \quad k_7=k_8=k_9=3, \quad k_{10}=k_{11}=4,$$
$$k_{12}=k_{13}=5, \quad k_{14}=k_{15}=6, \quad k_{50}=22, \quad k_{100}=45, \quad k_{1000}=485.$$
\par Finally, let us note that equality \eqref{nev_ksi_n_K_theta_n_ball_eq} holds
for those and only those dimensions $n$ when $k_n=1$.

\bigskip
\par\noindent{\bf Acknowledgements}
\bsk

I am very thankful to Pavel Shvartsman for useful suggestions and remarks.

\bigskip




\begin{thebibliography}{ABC}

\addtocontents{toc}{References \hfill\thepage\par}




\bibitem{barthelmann_2000} {\sc V. Barthelmann, E. Novak, and K. Ritter.} High dimensional polynomial interpolation on sparse grids, {\it Adv. Comput. Math.}, {\bf 12}:4 (2000), 273--288.
\bibitem{bauer_1958}
{\sc H. Bauer.} Minimalstellen von Funktionen und Extremal punkte,
 {\it Archiv der Mathematik}, {\bf 9}:4 (1958), 389--393.
\bibitem{deboor_1994}
{\sc C. de Boor.} {\it Polynomial Interpolation in Several Variables},  In: R. de Millo and J. R. Rice
(editors), Studies in Computer Science, 87--119, Plenum Press, 1994.
\bibitem{clements_lindstrem_1965} {\sc G.\,F. Clements and B. Lindstr\"om.} A sequence of $(\pm 1)$
determinants with large values, {\it Proc. Amer. Math. Soc.}, {\bf 16} (1965), 548--550.
\bibitem{demarchi_2015} {\sc S. De Marchi.} {\it  Lectures on Multivariate Polynomial Interpolation,}
G\"{o}ttingen\,--\,Padova, 2015.
\bibitem{devore_lorentz_1993} {\sc R.\,A. DeVore and G.\,G. Lorentz.}
{\it Constructive Approximation}, Springer-Verlag: Berlin\,--\,Heidelberg, 1993.
\bibitem{fejes_tot_1964} {\sc L. Fejes T\'{o}t.}
{\it Regular Figures}, New York: Macmillan/Pergamon, 1964.
\bibitem{fiht_2001} {\sc G.\,M. Fikhtengol'ts.}
{\it The Course in Differential and Integral Calculation. Vol.3,} Fizmatlit, Moscow, 2001 (in~Russian).
\bibitem{gasca_2000}
{\sc M. Gasca and T. Sauer.} Polynomial interpolation in several variables, {\it Adv. Comput. Math.},
{\bf 12}:4 (2000), 377--410.
\bibitem{gunzburger_2014} {\sc M. Gunzburger and A. Teckentrup.} Optimal point sets for total degree polynomial
interpolation in moderate dimensions, arXiv:1407.3291v1 [math.NA] 11 Jul 2014.
 \bibitem{horadam_2007} {\sc K.\,J. Horadam.} {\it Hadamard Matrices and Their Applications,} Princeton University Press, 2007.
\bibitem{hadamard_1893}
{\sc J.  Hadamard.}  R\'esolution d'une question relative aux d\'eterminants,
{\it Bull. Sciences Math. (2),} {\bf 17}\,(1893), 240--246.
\bibitem{hall_1970}
{\sc M. Hall, Jr.} {\it Combinatorial Theory}, Blaisdall Publishing Company,
 Waltham (Massa\-chusetts)\,--\,Toronto\,--\,London, 1967.
 \bibitem{hudelson_1996}
 {\sc M. Hudelson, V. Klee, and~D. Larman.} Largest $j$-simplices
in $d$-cubes: some relatives \linebreak of the~Hadamard maximum determinant
problem, {\it Linear Algebra Appl.},   {\bf 241--243} (2019), 519--598.
 \bibitem{kaarnioja_2015}
 {\sc V. Kaarnioja.} Multivariate polynomial interpolation over arbitrary grids,
 arXiv:1512.07424v1 [math.NA] 23 Dec 2015.
\bibitem{lassak_beitr_2011}
{\sc M.  Lassak.} Approximation of convex bodies by inscribed simplices of maximum volume,
{\it  Discr. Comput. Geom.,}  {\bf 21} (1999),  449--462.
\bibitem{manjhi_2022}
{\sc P.\,K. Manjhi and M.\,K. Rama.} Some new examples of circulant partial Hadamard matrices of type $4 - H(k\times n)$,
{\it Adv. Appl. Math. Sci.,} {\bf 21}:5\,(2022), 2559--2564.
\bibitem{nevskii_mais_2003_10_1}
{\sc M.\,V. Nevskii.} Estimates for the minimum norm of a projection in linear interpolation over the~vertices of an
$n$-dimensional cube, {\it Model. Anal. Inform. Sist.,} {\bf 10}:1 (2003), 9--19 \linebreak (in~Russian).
\bibitem{nevskii_mais_2006_13_2}
{\sc M.\,V. Nevskii.}  Geometric methods in the minimal projection problem,
 {\it Model. Anal. Inform. Sist.,} {\bf 13}:2 (2006), 16--29
(in~Russian).
\bibitem{nevskii_mais_2008_15_3}
{\sc M.\,V. Nevskii.}  Inequalities for the norms of interpolating projections,
 {\it Model. Anal. Inform. Sist.,} {\bf 15}:3 (2008), 28--37
(in~Russian).
\bibitem{nevskii_mais_2009_16_1}
{\sc M.\,V. Nevskii.}   On a certain relation for the minimal norm of an interpolational projection,
 {\it Model. Anal. Inform. Sist.,} {\bf 16}:1 (2009), 24--43
(in~Russian).
\bibitem{nevskii_matzam_2010}
M.\,V. Nevskii,  On a property of $n$-dimensional simplices, {\it Mat. Zametki,}
{\bf 87}:4 (2010), 580--593
(in~Russian). English translation: {\it  Math. Notes,} {\bf 87}:4 (2010), 543--555.
\bibitem{nevskii_mais_2011_2}
{\sc M.\,V. Nevskii.} On geometric characteristics of an $n$-dimensional simplex,
 {\it Model. Anal. \linebreak  Inform. Sist.,} {\bf 18}:2 (2011), 52--64
(in~Russian).
\bibitem{nevskii_dcg_2011}
{\sc M. Nevskii.} Properties of axial diameters of a simplex,
{\it  Discr. Comput. Geom.,}  {\bf 46}:2 (2011),  301--312.
\bibitem{nevskii_monograph}
{\sc M.\,V. Nevskii.}
{\it Geometric Estimates in Polynomial
Interpolation}, P.\,G.~De\-mi\-dov Yaroslavl'  State University, Yaroslavl', 2012 (in~Russian).
\bibitem{nevskii_mais_2018_25_6}
{\sc M.\,V. Nevskii.}~On some problems for a simplex and a ball in ${\mathbb R}^n$,
{\it Model. Anal. Inform. Sist.}, {\bf 25}:6 (2018), 680--691
(in~Russian).
English translation: {\it Aut.
Control Comp. Sci.}, {\bf 53}:7 (2019), 644--652.
\bibitem{nevskii_mais_2019_26_3}
{\sc M.\,V. Nevskii.} Geometric estimates in interpolation on an n-dimensional ball,
{\it Model. Anal. Inform. Sist.}, {\bf 26}:3 (2019), 441--449
(in~Russian).
 English translation: {\it Aut.
Control Comp. Sci.},  {\bf 54}:7 (2020), 712--718.
\bibitem{nevskii_mais_2021_28_2}
{\sc M.\,V. Nevskii.} On properties of a regular simplex inscribed into a ball,
{\it Model. Anal. Inform. Sist.}, {\bf 28}:2 (2021), 186--197
(in~Russian).
 English translation: {\it Aut.
Control Comp. Sci.},  {\bf 56}:7 (2022), 778--787.
\bibitem{nevskii_mais_22}
{\sc M.\,V. Nevskij.} On some estimate for the norm of an interpolation projector, {\it  Model. Anal. Inform. Sist.,}{\bf 29}:2 (2022), 92--103 (in~Russian).  English translation: {\it  Aut. Control Comp. Sci.,} {\bf 57}:7 (2023), 718--726.
\bibitem{nevskii_matzam_23}
{\sc  M.\,V. Nevskii.} On the minimal norm of the projection operator for linear interpolation
\linebreak on~an~$n$-dimensional ball, {\it  Mat. Zametki,} {\bf114}:3 (2023), 477--480
(in~Russian). English translation: {\it Math. Notes,} {\bf 114}:3 (2023), 415--418.
\bibitem{nevskii_arxiv_2024}
{\sc M. Nevskii.} Geometric estimates in linear interpolation on a cube and a ball,
arXiv:2402.11611v1 [math.CA] 18 Feb 2024.


\bibitem{nev_ukh_mais_2016_23_5}
{\sc M.\,V. Nevskii  and A.\,Yu. Ukhalov.} On~numerical characteristics of a simplex and their estimates,
{\it Model. Anal. Inform. Sist.}, {\bf 23}:5 (2016), 602--618
(in~Russian).
English translation: {\it Aut.
Control Comp. Sci.}, {\bf 51}:7 (2017), 757--769.
\bibitem{nev_ukh_beitrage_2018}
{\sc M. Nevskii and A. Ukhalov.}  Perfect Simplices in ${\mathbb R}^5$,
{\it  Beitr. Algebra Geom.,}  {\bf 59}:3 (2018), 501--521.
\bibitem{nev_ukh_mais_2018_25_3}
{\sc M.\,V. Nevskii  and A.\,Yu. Ukhalov.}
On optimal interpolation by linear functions on an $n$-dimensional cube,
{\it Model. Anal. Inform. Sist.}, {\bf 25}:3 (2018), 291--311
(in~Russian).
English translation: {\it Aut.
Control Comp. Sci.},  {\bf 52}:7 (2018), 828--842.
\bibitem{nev_ukh_mais_2019_26_2}
{\sc M.\,V. Nevskii  and A.\,Yu. Ukhalov.}
 Linear interpolation on a Euclidean ball in ${\mathbb R}^n$,
{\it Model. Anal. Inform. Sist.}, {\bf 26}:2 (2019), 279--296
(in~Russian).
English translation: {\it Aut.
Control Comp. Sci.},  {\bf 54}:7 (2020), 601--614.
\bibitem{nev_ukh_posobie_2020}
{\sc M.\,V. Nevskii  and A.\,Yu. Ukhalov.}
{\it Selected Problems in Analysis and Computational \linebreak Geometry. Part 1},
  P.\,G.~De\-mi\-dov Yaroslavl'  State University, Yaroslavl', 2020 (in~Russian).
 \bibitem{nev_ukh_posobie_2022}
{\sc M.\,V. Nevskii  and A.\,Yu. Ukhalov.}
{\it Selected Problems in Analysis and Computational \linebreak Geometry. Part 2},
  P.\,G.~De\-mi\-dov Yaroslavl'  State University, Yaroslavl', 2022 (in~Russian).
\bibitem{pashkovskij_1983}
{\sc S. Pashkovskij.}
{\it Vychislitel'nye Primeneniya Mnogochlenov i Ryadov Chebysheva}, Nauka, Moskva, 1983 (in~Russian).
\bibitem{prudnikov_1981}
{\sc A.\,P. Prudnikov, Yu.\,A. Brychkov, and O.\,I. Marichev.}
{\it Integraly i Ryady,} Nauka, Moskva, 1981
(in~Russian).
\bibitem{rivlin_1974}
{\sc T.\,J. Rivlin.} {\it The Chebyshev Polynomials,} John Wiley \& Sons,
New York\,--\,Chichester\,--\,Brisbane\,--\,Toronto, 1974.
\bibitem{scott_1985}
{\sc  P.\,R. Scott.}
Lattices and convex sets in space, {\it Quart. J. Math. Oxford Ser. (2),} {\bf 36} (1985), 359--362.
\bibitem{scott_1989}
 {\sc P.\,R. Scott.}
Properties of axial diameters, {\it Bull. Austral. Math. Soc.,} {\bf 39} (1989), 329--333.
\bibitem{slepian_1969}
{\sc D. Slepian.}
The content of some extreme simplices,
{\it Pacific J.~Math.,} {\bf 31}\,(1969), 795--808.
\bibitem{suetin_1979}
{\sc P.\,K. Suetin.} {\it Klassicheskie ortogonal'nye mnogochleny}, Nauka, Moskva,1979 (in~Russian).
\bibitem{sege_1975} {\sc G. Szeg\"{o}.} {\it Orthogonal Polynomials}, Providence, American Mathematical Society,  New~York, 1975.
\bibitem{vandev_1992}
{\sc D.  Vandev.}
A minimal volume ellipsoid around a simplex,
{\it C. R. Acad. Bulg.~Sci.}, {\bf 45}:6\,(1992), 37--40.
\end{thebibliography}
\end{document}